\documentclass[11pt]{article}
\usepackage{epsfig,graphics,empheq, graphicx,amsmath,amssymb,amsfonts}
\usepackage{epstopdf}
\usepackage{color}
\usepackage{amsmath}
\usepackage{cases}
\usepackage[position=bottom]{subfig}
\usepackage{paralist}
\usepackage{enumitem}
\usepackage{graphicx}
\usepackage[percent]{overpic}
\usepackage{tikz}

\newsavebox{\largestimage}

\setlist{nolistsep}

\newcommand{\bm}[1]{\mbox{\boldmath{$#1$}}}

\def\u{{\bm u}}
\def\v{{\bm{v}}}
\def\x{{\bm{x}}}

\def\0{{\bm 0}}

\def\cl {\nonumber \\}
\def\el {\nonumber}

\DeclareMathOperator*{\argmin}{arg\,min}

\begin{document}

\title{A kinetic theory approach to model crowd dynamics with disease contagion}

\author{Daewa Kim$^*$ and Annalisa Quaini$^{**}$ \\
\footnotesize{
$^*$Department of Mathematics, West Virginia University, 94 Beechurst Ave, Morgantown, WV 26506}\\
\footnotesize{
$^{**}$Department of Mathematics, University of Houston, 3551 Cullen Blvd, Houston TX 77204}\\
\footnotesize{daewa.kim@mail.wvu.edu; quaini@math.uh.edu}
}

\maketitle

\noindent{\bf Abstract}
We present some ideas on how to extend a kinetic type model for crowd dynamics
to account for an infectious disease spreading. We focus on a medium size
crowd occupying a confined environment where the disease
is easily spread. The kinetic theory approach we choose uses 
tools of game theory to model the interactions of a person with the surrounding people 
and the environment and it features a parameter to represent the level of stress. 
It is known that people choose different walking strategies when subjected to fear
or stressful situations. 
To demonstrate that our model for crowd dynamics could be used to reproduce realistic scenarios, 
we simulate passengers in one terminal of Hobby Airport in Houston.
In order to model disease spreading in a walking crowd, we introduce a variable that denotes 
the level of exposure to people spreading the disease. In addition, we introduce 
a parameter that describes the contagion interaction strength and 
a kernel function that is a decreasing function of the distance between a person and a spreading individual.
We test our contagion model on a problem involving a small crowd walking through a corridor.

%% section 1 %%%%%%%%%%%%%%%%%
\section{Introduction} \label{sec:Intro}
%%%%%%%%%%%%%%%%%%%%%%%%

We are interested in studying the early stage of an infectious disease spreading in an intermediate size
population occupying a confined environment, such as an airport terminal or a school, 
for a short period of time (minutes or hours). 
Classical models in epidemiology use mean-field approximations
based on averaged large population behaviors over a long time span (typically weeks or months).
Obviously, such models fail when population size is small-to-medium. 
Our overarching goal is to model the spreading of a disease in a walking crowd by
extending a kinetic theory approach for crowd dynamics that compares favorably with 
experimental data for a medium-sized population~\cite{kim_quaini}.
In this paper, we present the key features of our crowd dynamics model
and preliminary ideas for the extension, that are tested in 1D cases.
We assume that the disease is such that it spreads with close proximity of individuals, like, e.g. measles, 
influenza or COVID-19.

The reason why we focus on a mesoscopic model for crowd dynamics is related to 
our interest in confined environments and medium-sized crowds. In broad terms, the large variety of models proposed over the years
can be divided into three main categories depending
on the (microscopic, mesoscopic, or macroscopic) scale of observation \cite{Bellomo2011}. 
Macroscopic models (see, e.g., \cite{Cristiani2014_book,Hughes2003,5773492}) treat the crowd as a continuum flow, 
which is well suited for large-scale, dense crowds. Microscopic models (see, e.g.,
\cite{ASANO2010842,bandini2009,Chraibi2011425,Chraibi2019,DAI20132202,Helbing1995, 6701214,6248013} and references therein)
use Newtonian mechanics to interpret pedestrian movement as the physical interaction between the people and the environment.
Mesoscale models (see, e.g.,~\cite{Agnelli2015,Bellomo2011383,Bellomo2017_book,Bellomo2013_new,Bellomo2015_new,Bellomo2016_new,Bellomo2019_new,Bellomo2013,Bellomo2012}) use a Boltzmann-type evolution equation for the statistical
distribution function of the position and velocity of pedestrians,
in a framework close to that of the kinetic theory of gases.
There is one key difference though: the interactions in Boltzmann equations 
are conservative and reversible, while the interactions in the kinetic theory of active particles are irreversible, 
non-conservative and, in some cases, non-local and nonlinearly additive. 
An important consequence is that often for active particles the
Maxwellian equilibrium does not exists~\cite{Aristov_2019}.
Another reason why we choose to work with a kinetic type model is the flexibility
in accounting for multiple interactions (hard to achieve in microscopic models) 
and heterogeneous behavior in people (hard to achieve in macroscopic models).
Finally, we would like to mention that multiscale approaches are possible as well. 
See, e.g., Ref.~\cite{Bellomo2020} for a multiscale vision to human crowds which provides a consistent 
description at the three possible modeling scales.

The first paper of this manuscript is dedicated to a description of a crowd dynamics model. The model, first presented 
in \cite{kim_quaini}, is based on earlier works \cite{Agnelli2015,Bellomo2015_new}
and is capable of handling evacuation from a room with one or more exits
of variable size and in presence of obstacles. 
The main ingredients of the model are the following: (i) discrete velocity directions
to take into account the granular feature of crowd dynamics, (ii) interactions modeled using
tools of stochastic games, and (iii) heuristic, deterministic modeling of the speed corroborated by
experimental evidence~\cite{Schadschneider2011545}. In \cite{kim_quaini}, we show that 
for groups of 40 to 138 people the average people density and flow rate computed
with our kinetic model is in great agreement with the respective measured quantities 
reported in a recent empirical study focused on egressing from a facility \cite{Kemloh}. 

To demonstrate that our model for crowd dynamics could be used to reproduce realistic scenarios, 
we simulate passengers in one terminal of Hobby Airport in Houston (USA). In a first set of tests, 
the passengers from two planes at the opposite ends of the terminal walk through the terminal
to reach the exit. In the second set of tests, we add a group of passengers that enters the terminal
through the entrance at the same time as the other two groups deplane and is directed to a gate. 
The aim of both sets of tests is to understand how the presence of obstacles in the terminal 
affects the egress time. Obviously, the longer one stays in the terminal in close proximity
with other individuals the more likely he/she gets infected. Thus, the egress time
from a potentially crowded confined environment, such as an airport terminal, 
is a key factor in the early spreading of an airborne disease. This is why we chose
these tests and how they are connected to the second part of the paper. 

In order to model disease spreading in a walking crowd, 
we take inspiration from from the work on emotional contagion (i.e., spreading of fear or panic) 
in Ref.~\cite{Bertozzi2015}.
We introduce a variable that denotes the level of exposure to people spreading the disease, with the underlying
idea that the more a person is exposed the more likely they are to get infected.
The model features a parameter that describes the contagion
interaction strength and a kernel function that is a decreasing function of the distance between
a person and a spreading individual. As a simplification, we assume that walking speed
and direction are given. We will show preliminary results for a problem involving a small crowd walking through a corridor.
The simplifying assumption will be removed in a follow-up paper, 
where the people dynamics will be provided by the complex pedestrian model described 
in the first part of the paper.
The approach we have in mind is different from what we used in Ref.~\cite{kim_quaini2020}. Therein
we coupled the pedestrian dynamics model to a disease contagion model, while in the future
we intend to add to the pedestrian dynamics model terms that account for disease spreading.

For related work on coupled dynamics of virus infection and healthy cells
see, e.g. Ref.~\cite{doi:10.1137/19M1250261} and references therein. 
A multiscale model of virus pandemic accounting for the interaction of different spatial scales (from the
small scale of the virus itself and cells, to the large scale of individuals and further up to
the collective behavior of populations) is presented in Ref.~\cite{Bellomo_preprint}.

The paper is organized as follows. 
Sec.~\ref{sec:ped_model} describes the crowd dynamics model, its full discretization, 
and shows numerical results in an airport terminal. 
In Sec.~\ref{sec:contagion}, we introduce our simplified contagion model.  
The discretization and preliminary results are also shown in Sec.~\ref{sec:contagion}.
Conclusions are drawn in Sec.~\ref{sec:concl}.

\section{A kinetic model for crowd dynamics}\label{sec:ped_model}
%%%%%%%%%%%%%%%%%

%The model we consider is based on the model proposed in \cite{Agnelli2015}.
Let  $\Omega \subset \mathbb{R}^2$ denote a bounded domain where people 
are walking to reach an exit $E$ that is either within the domain or belongs 
to the boundary $\partial \Omega$.
The case of multiple exits (i.e., $E$ is the finite union of disjoint sets) can be easily handled as well. 
The rest of the boundary is made of walls, denoted with $W$. Walls and other kinds of obstacles  
could be present within the domain. 
Let $\x=(x,y)$ denote position and ${\v}=v (\cos \theta, \sin \theta) $ denote  
velocity, where $v$ is the velocity modulus and $\theta$ is the velocity direction.
 For a large group of people inside $\Omega$, let 
 \[ f= f(t, \x, v, \theta)\quad \text{for all} \,\,\, t \ge 0,  \,\, \x \in \Omega, \,\, v \in [0, V_M], \,\, \theta \in [0, 2\pi), \]
 where $V_M$ is the largest speed a person can reach in low density and optimal environmental conditions.
Under suitable integrability conditions, $f(t, \x, v, \theta)d \x d v d \theta$ represents the number of individuals who, at time $t$, 
are located in the infinitesimal rectangle $[x, x+dx] \times [y, y+dy]$ and have a velocity belonging to $[v, v + dv] 
\times [\theta, \theta+d\theta]$. 

For simplicity and following \cite{Agnelli2015}, we make two simplifying assumptions
on the velocity vector:
\begin{enumerate}
\item Variable $\theta$ is discrete, i.e.~it can take values in the set:
\[ I_{\theta}=\left \{ \theta_{i}= \frac{i-1}{N_d} 2\pi : i = 1, \dots, N_d \right \}, \]
where $N_d$ is the maxim number of possible directions. 
\item People adjust their walking speed $v$ depending on the level of congestion around them, i.e.~$v$ is treated as a
deterministic variable.
\end{enumerate}
The second assumption is corroborated by experimental studies that show that
the walking speed mainly depends on the local level of congestion.
Given the deterministic nature of the variable $v$, the distribution function can be written as
\begin{equation}\label{eq:f}
f(t, \x, \theta)= \sum_{i=1}^{N_d} f^i(t, \x)\delta(\theta - \theta_i),
\end{equation}
where $f^i(t, \x)=f(t, \x, \theta_i)$ represents the people that, at time $t$  and position $\x$, move 
with direction $\theta_i$. In eq.~\eqref{eq:f}, $\delta$ denotes the Dirac delta function.

In the rest of the paper, we will work with dimensionless variables. To this purpose, 
we introduce the following reference quantities: 
\begin{itemize}
\item[-] $D$: the largest distance a pedestrian can cover in domain $\Omega$;
\item[-] $T$: a reference time given by $D/V_M$ (recall that $V_M$ is the largest speed a person can reach),
\item[-] $\rho_M$: the maximal admissible number of pedestrians per unit area.
\end{itemize}
The dimensionless variables are then: position $\hat{\x}=\x/D$, time $\hat{t}=t/T$, velocity modulus $\hat{v}=v/V_M$ and distribution function $\hat{f}=f/ \rho_M$. For ease of notation, the hats will be omitted with the understanding that all variables are
dimensionless from now on. 

Due to the normalization of $f$ ant the $f_i$, $i = 1, \dots, N_d$, the dimensionless local density is obtained 
by summing the distribution functions over the set of directions:
\begin{align}\label{eq:rho}
\rho(t, \x)=\sum_{i=1}^{N_d}f^i(t, \x) .
\end{align}
Following assumption 2 mentioned above, the walking speed is given by $v=v[\rho](t, \x)$, 
where square brackets are used to denote that $v$ depends on $\rho$ in a functional way. 
For instance, $v$ can depend on $\rho$ and on its gradient. 

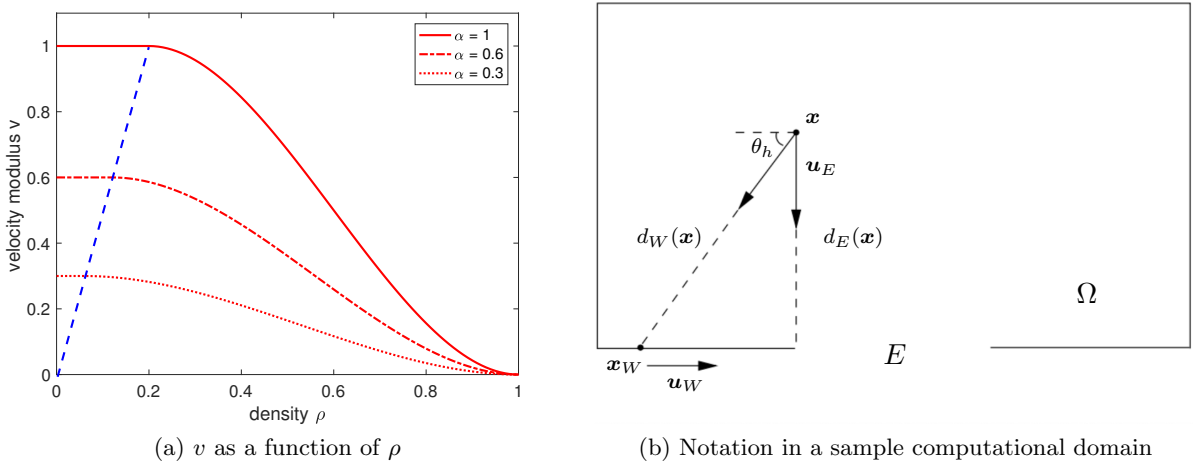
\begin{figure}[h!]
\centering
\subfloat[$v$ as a function of $\rho$]{
\begin{overpic}[width=0.48\textwidth,grid=false]{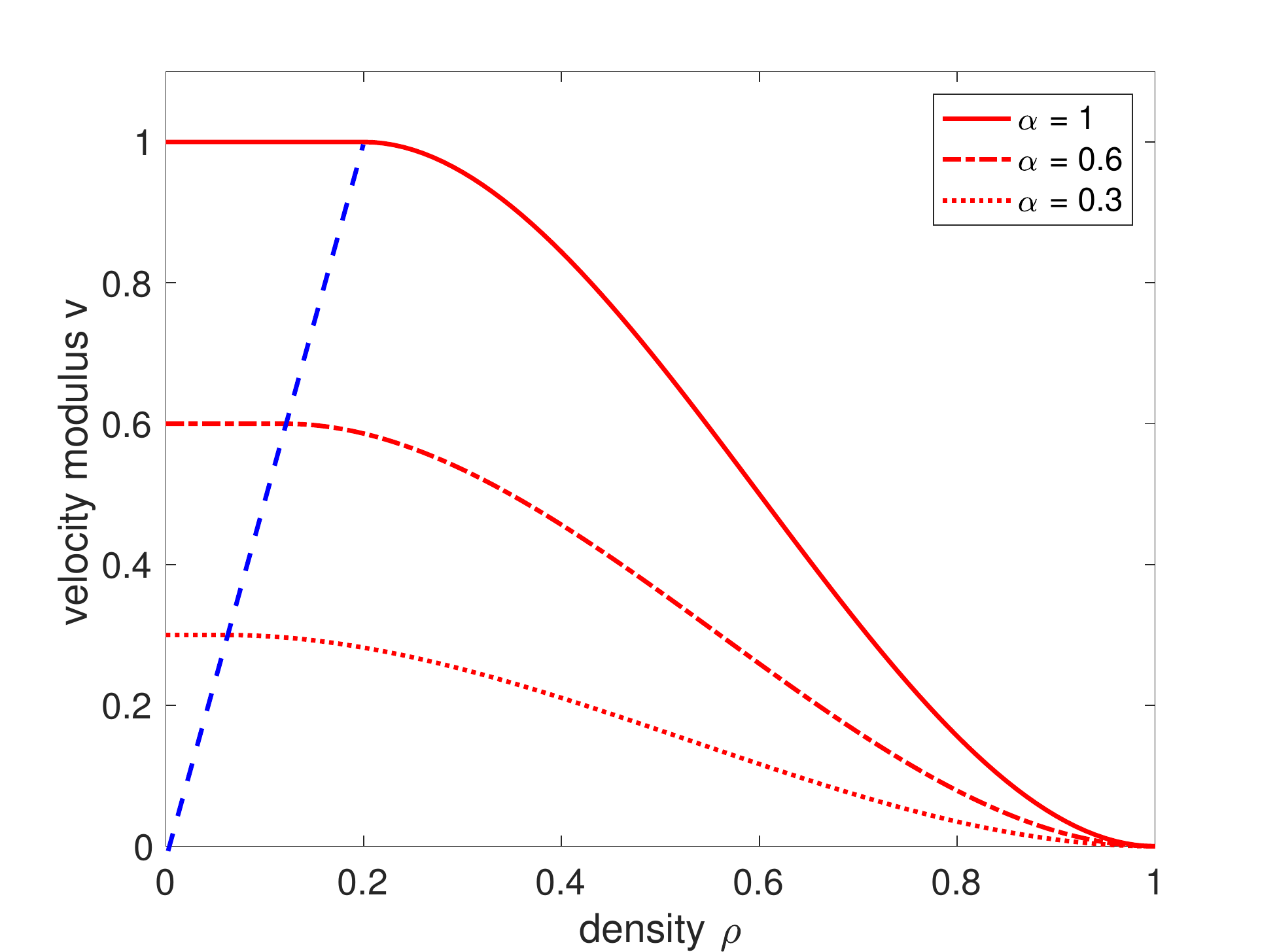}
\end{overpic}
}
\subfloat[Notation in a  sample computational domain]{
\begin{overpic}[width=0.485\textwidth, grid=false, tics=10]{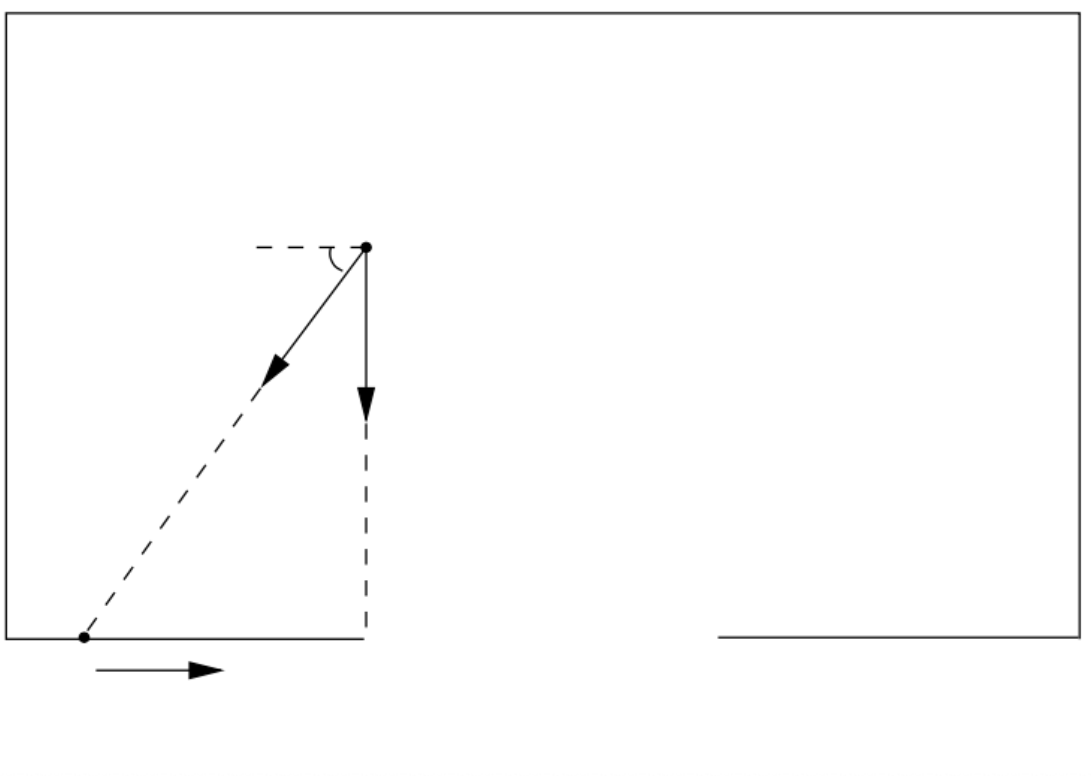}
      \put(80,20){\text{$\Omega$}}
      \put(48,10){\text{$E$}}
      \put(35,50){\scriptsize{\text{$\x$}}}
      \put(2, 9){\scriptsize{\text{$\x_{W}$}}}
      \put(38,30){\scriptsize{\text{$d_E(\x)$}}}
      \put(7,30){\scriptsize{\text{$d_W(\x)$}}}
      \put(35,42){\scriptsize{\text{$\u_E$}}}
      \put(12,6){\scriptsize{\text{$\u_W$}}}
     \put(26,45){\scriptsize{\text{$\theta_h$}}}
\end{overpic}
}
\caption{(a) Dependence of the dimensionless walking speed $v$ on the dimensionless density $\rho$ 
for different values of the parameter $\alpha$, which represents the quality of the walkable domain.
(b) Sketch of computational domain $\Omega$ with exit $E$ and a pedestrian located at $\x$, moving
with direction $\theta_h$. The pedestrian should choose direction $\u_E$ to reach the exit, 
while direction $\u_W$ is to avoid collision with the wall. The distances form the exit and from the wall
are $d_E$ and $d_W$, respectively.}\label{velocity}
\end{figure}

In order to define the walking speed, we introduce a parameter $\alpha \in [0, 1]$ to represent the quality of the 
walkable domain: where $\alpha=1$ people can walk at the desired speed (i.e., $V_M$) because the domain is clear, 
while where $\alpha=0$ people are forced to slow down or stop because an obstruction is present. 
For simplicity, we assume that the maximum dimensionless speed a person can reach
is equal to $\alpha$.
Let $\rho_c$ be a critical density 
value such that for $\rho < \rho_c$ we have free flow regime (i.e.,~low density condition), while for $\rho > \rho_c$
we have a slowdown zone (i.e.,~high density condition). Following the experimentally measured values of $\rho_c$ reported 
in \cite{Schadschneider2011545}, we set $\rho_c = \alpha/5$.   
Then, we set the walking speed $v$ equal to $\alpha$ in the free flow regime
and equal to a heuristic third-order polynomial in the slowdown zone:
\begin{align}\label{eq:v}
v=v(\rho)=
\begin{cases}
\alpha \quad & \text{for} \quad \rho \leq \rho_c(\alpha)= \alpha/5  \\
a_3\rho^3+a_2\rho^2+a_1\rho+a_0 \quad & \text{for} \quad \rho > \rho_c(\alpha)=\alpha/5   ,     
\end{cases}
\end{align}
where $a_i$ is constant for $i = 0,1,2,3$.
To set the value of these constants, we impose the following 
conditions: $v(\rho_c) = \alpha$,  $\partial_{\rho} v(\rho_c) = 0$, $v(1)=0$ and $\partial_{\rho} v(1) = 0$, which 
lead to:
\begin{align}\label{eq:coeff}
\begin{cases}
a_0 &= (1/(\alpha^3-15\alpha^2+75\alpha-125))(75\alpha^2-125\alpha) \\
a_1 &= (1/(\alpha^3-15\alpha^2+75\alpha-125))(-150\alpha^2)\\
a_2 &= (1/(\alpha^3-15\alpha^2+75\alpha-125))(75\alpha^2+375\alpha) \\
a_3 &= (1/(\alpha^3-15\alpha^2+75\alpha-125))(-250\alpha). 
\end{cases}
\end{align}
%\daewa{Daewa, make sure to update the figure for copyright reasons.} 
Fig.~\ref{velocity} (a) shows $v$ as a function of $\rho$ for $\alpha = 0.3, 0.6, 1$.

%%%%%%%%%%%%%%%%%%
%\begin{figure}[h!tb]
%\centering
%\begin{overpic}[width=0.5\textwidth,grid=false]{./velocity_new_0424.pdf}
%\end{overpic}
%\caption{Dependence of the dimensionless velocity modulus $v$ on the dimensionless density $\rho$ for different values of the parameter $\alpha$, which represents the quality of the environment.}
%\label{velocity}
%\end{figure}

%%%%%%%%%%%%%%%%%%
\subsection{Modeling interactions} \label{modelinginteractions}

Let us consider the scenario depicted in Fig.~\ref{velocity} (b), where there is a person
located a point $\x$ that needs to reach exit $E$. We model the path 
this person takes as the result of four factors:
\begin{enumerate}[label={(F\arabic*})]
\item The goal to reach the exit.
\item The desire to avoid collisions with the walls.
\item The tendency to look for less congested area.
\item The tendency to follow the stream or herding.
\end{enumerate}
Factors F1 and F2 are related to geometric aspects of the domain, while factors F3 and F4
consider that people's behavior is strongly affected by the surrounding crowd. 
These last two factors are dominant in different situations: F4 emerges in stressful situations, 
while F3 characterizes normal behavior. 
To weight between F3 and F4, we use parameter $\varepsilon \in [0,1]$ with $\varepsilon=0$ 
(resp., $\varepsilon=1$) if F3 (resp., F4) prevails.

In order to explain how the four factors are modeled, we need to introduce some
terminology. Interactions involve three types of people: 
\begin{itemize}
\item[-] \textit{test people} with state $(\x, \theta_i)$: they are representative of the whole system;
\item[-] \textit{candidate people} with state $(\x, \theta_h)$: they can reach in probability the state of the test people 
after individual-based interactions with the environment or with field people; 
\item[-] \textit{field people} with state $(\x, \theta_k)$: their interactions with candidate people triggers a possible change of state.
\end{itemize}
We note that while the candidate person modifies their state, in probability, into that of the test person due to interactions with field people, 
the test person loses their state as a result of these interactions.

Next, we introduce some notation.
Given a candidate person at point $\x$ in the walkable domain $\Omega$, we define its distance to the exit as
\[d_E(\x)= \min_{\x_E \in E} || \x-\x_E ||,\]
and we consider the unit vector $\u_E(\x)$, pointing from $\x$ to the exit. See Fig.~\ref{velocity} (b).
Both $d_E$ and $\u_E$ will be used to model (F1). 

Assume that the candidate person at $\x$ is moving with direction $\theta_h$. 
We define the distance $d_W(\x, \theta_h)$ from the person to a wall at a point $\x_W(\x, \theta_h)$
where the person is expected to collide with the wall if they do not change direction.
The unit tangent vector $\u_W(\x, \theta_h)$ to $\partial \Omega$ at $\x_W$
points to the direction of the exit. See Fig.~\ref{velocity} (b).
Vector $\u_W$ is used to avoid a collision with the walls, i.e.~to model (F2). 
%\anna{Daewa, plaese add $\x_W$ to the figure. Also adjust Fig.~\ref{velocity} (b) so that the two . are aligned.}

In order to model (F3), i.e.~the decision of candidate person $(\x, \theta_h)$ to change direction in order to avoid congested areas, 
we use the direction that gives the minimal directional derivative of the density at the point $\x$. 
We denote such direction by unit vector $\u_C(\theta_h, \rho)$.

Finally, we introduce unit vector $\u_F=(\cos\theta_k, \sin\theta_k)$ to model (F4), i.e.~the decision of candidate person with direction $\theta_h$ to follow a field person with direction $\theta_k$
with whom they came into contact.

%%%%%%%%%%%%%%%%%%%%%%
%\begin{figure}[h!]
%\centering
%\begin{overpic}[width=0.6\textwidth,grid=false,tics=10]{./domain_0409}
%\end{overpic}
%\caption{A particle located at point $\x$.}
%\label{domain}
%\end{figure}
%%%%%%%%%%%%%%%%%%%%%%%

%%%%%%%%%%%%%%%%%%%%%%%%
\subsubsection{Interaction with the walls}

We assume that people change direction, in probability, only to an adjacent 
clockwise or counterclockwise direction in set $I_\theta$. This means a candidate person 
with walking direction $\theta_h$ may choose directions $\theta_{h-1}, \theta_{h+1}$ 
or keep direction $\theta_h$. For $h=1$ we set $\theta_{h-1}=\theta_{N_d}$ and 
for $h=N_d$ we set $\theta_{h+1}=\theta_1$.
Let $\mathcal{A}_h(i)$ be the \textit{transition probability}, i.e.~the probability that a candidate person with direction $\theta_h$ 
adjusts their direction to $\theta_i$ (the direction of the test person) due to the presence of walls and/or an exit. 
The following constraint for $\mathcal{A}_h(i)$ has to be satisfied:
\[
\sum_{i=1}^{N_d} \mathcal{A}_h(i)=1 \quad \text{for all} \,\, h \in \{1, \dots, N_d\}.
\]
The set of all transition probabilities $\mathcal{A}=\{\mathcal{A}_h(i) \}_{h,i= 1, \dots, N_d}$ 
forms the so-called \textit{table of games} that models the game played by active people interacting with the walls.

We define the vector
\begin{equation}\label{eq:uG}
\u_G(\x, \theta_h)= \frac{(1-d_E(\x))\u_E(\x) + (1-d_W(\x, \theta_h))\u_W(\x, \theta_h)}{|| (1-d_E(\x))\u_E(\x) + (1-d_W(\x, \theta_h))\u_W(\x, \theta_h) ||} = (\cos \theta_G, \sin \theta_G).
\end{equation}
Here $\theta_G$ is the  \textit{geometrical preferred direction}, which is the ideal direction 
that a person should take in order to reach the exit (factor F1) and avoid the walls (factor F2)
in an optimal way. Notice that the closer a person is
to an exit (resp., a wall), the more direction $\u_E$ (resp., $\u_W$) weights.

A candidate person with direction $\theta_h$ will change their direction by choosing
the angle closest to $\theta_G$ among the three allowed 
directions $\theta_{h-1}, \theta_{h}$ and $\theta_{h+1}$. The transition probability is given by:
\begin{equation}\label{eq:A}
\mathcal{A}_h(i)=\beta_h(\alpha)\delta_{s,i} + (1-\beta_h(\alpha))\delta_{h, i}, \quad i=h-1, h, h+1,
\end{equation}
where
\[s=\argmin_{j \in \{h-1,h+1\}}\{d(\theta_G, \theta_j)\},\]
with
\begin{equation}\label{eq:distance}
 d(\theta_p, \theta_q)=
\begin{cases}
|\theta_p - \theta_q|  & \text{if} \,\,\, |\theta_p - \theta_q | \leq \pi, \\
2\pi - |\theta_p- \theta_q|  & \text{if} \,\,\, |\theta_p- \theta_q| > \pi .             
\end{cases}
\end{equation}
In \eqref{eq:A}, $\delta$ denotes the Kronecker delta function. Coefficient $\beta_h$ is defined by:
\[
\beta_h(\alpha)=
\begin{cases}
\alpha & \text{if} \,\,\, d(\theta_h, \theta_G) \geq \Delta\theta, \\
\alpha \dfrac{d(\theta_h, \theta_G)}{\Delta\theta}& \text{if} \,\,\, d(\theta_h, \theta_G)< \Delta \theta ,             
\end{cases}
 \]
where $\Delta\theta=2\pi/{N_d}$. The role of $\beta_h$ is to allow for a transition to $\theta_{h-1}$ or
$\theta_{h+1}$ even in the case that the 
geometrical preferred direction $\theta_G$ is closer to $\theta_h$. Such a transition is more likely 
to occur the more distant $\theta_h$ and $\theta_G$ are.
%Notice that if $\theta_G=\theta_h$, then $\beta_h=0$ and $\mathcal{A}_h(h)=1$, meaning that a 
%pedestrian keeps the same direction (in the absence of interactions other than with the 
%environment) with probability 1.

%%%%%%%%%%%%%%%%%%%%%%%%
\subsubsection{Interaction with obstacles}\label{sec:obstacles}

%The strategy reported in the previous section to avoid collisions with the walls works well
%when the pedestrian is sufficiently far from the walls. If pedestrians get too close to the 
%bounding walls, and in particular if they are close to an exit, the definition of $\u_G$ in \eqref{eq:uG}
%does not prevent collisions with the walls. Thus, obstacles within the domain $\Omega$
%cannot be avoided just by adjusting $\u_W$. In this section, we report
%an effective strategy to handle obstacles. 

In \cite{kim_quaini}, we introduced a strategy to handle obstacles 
within domain $\Omega$. This strategy uses three ingredients
to exclude the real obstacle area from the walkable domain:
\begin{enumerate}
\item An effective area: an enlarged area that encloses the real obstacle.
\item A definition of $\u_W$ to account for the effective area. 
\item A setting of the parameter $\alpha$ in the effective area depending on the shape of the obstacle. 
\end{enumerate}
The effective area is necessary especially if the obstacle is close to an exit: it allows to 
define $\u_W$ with respect to a larger area than the obstacle area itself to 
achieve the goal of having no people walk on the real obstacle area. 
In \cite{kim_quaini}, we found that the goal is successfully achieved 
with an effective area four times bigger than the real obstacle area.

Since some pedestrians will walk on part of the effective area, one needs to set parameter $\alpha$
in a suitable way. For a discussion oh how to set $\alpha$ to realize different obstacle shapes, 
we refer to \cite{kim_quaini}.

%By setting $\alpha=1$ (i.e., best quality of the environment) in the effective area, 
%pedestrians can move with the maximal velocity modulus as they approach the obstacle
%and thus they quickly adapt to the effective area through $\u_W$. However, some pedestrians
%will walk close to the top, bottom, and rear (with respect to the pedestrian motion) 
%boundary of the effective area. Thus, the real obstacle is located at the front of the effective area.
%From the numerical results reported in Sec.~\ref{sec:obstacle}, we also see that the shape
%of the obstacle is square. By setting $\alpha=0$ (i.e., worst quality of the environment) in the effective area, 
%pedestrians are forced to slow down at the front part of the effective area. The slow down leads
%to higher densities in the front part of the effective area, therefore direction $\u_W$ competes
%with direction $\u_C$. As a result some pedestrians walk on the front part of the effective area.
%However, as the congestion decreases pedestrians avoid the rear part of the effective area.
%From the numerical results shown in Sec.~\ref{sec:obstacle}, we see that the shape
%of the obstacle for $\alpha=0$ in the effective area is slender. 

%%%%%%%%%%%%%%%%%%%%%%%%%%%%% 
\subsubsection{Interactions between pedestrians}

As a candidate person with direction $\theta_h$ walks, they interact with a field person 
that moves with direction $\theta_k$. As a result of this interaction, the candidate 
person can change their direction to $\theta_i$ (direction of the test person) in the 
search for less congested areas if their stress level is low or change to $\theta_k$ (direction of the field person)
if their stress level is high.
The \textit{transition probability} is given by $\mathcal{B}_{hk}(i)[\rho]$.
The following constrain for $\mathcal{B}_{hk}(i)$ has to be satisfied:
\[
\sum_{i=1}^{N_d} \mathcal{B}_{hk}(i)[\rho]=1 \quad \text{for all} \,\, h, k  \in \{1, \dots, N_d\},
\]
where again the square brackets denote the dependence on the density $\rho$.
Of course, we are still under the assumption that people change direction, in probability, only to an adjacent 
clockwise or counterclockwise direction in set $I_\theta$.

To take into account the search for a less congested area (factor F3) and the tendency to herd (factor F4), 
for a candidate person with direction $\theta_h$ interacting with a field person with direction 
$\theta_k$ we define the vector
\[\u_P(\theta_h, \theta_k, \rho)= \frac{\varepsilon\u_F+(1-\varepsilon)\u_C(\theta_h, \rho)}
{||\varepsilon\u_F+(1-\varepsilon)\u_C(\theta_h, \rho)||} = (\cos \theta_P, \sin \theta_P),
\]
where the subscript $P$ stands for \textit{people}. Direction $\theta_P$ is the \textit{interaction-based preferred direction}, 
obtained as a weighted combination between the direction of the field person (i.e., $\u_F = (\cos\theta_k, \sin\theta_k)$)
the direction pointing to a less crowded area (i.e., $\u_C$). The latter direction can be computed 
for a candidate pedestrian with direction $\theta_h$ and located at $\x$, by taking 
\[C=\argmin_{j \in \{h-1, h, h+1\}}\{\partial_j\rho(t, \x)\},\]
where $\partial_j\rho$ denotes the directional derivative of $\rho$ in the direction given by angle $\theta_j$. 
We have $\u_C(\theta_h, \rho)=(\cos\theta_C, \sin\theta_C)$.

The transition probability is given by:
\[\mathcal{B}_{hk}(i)[\rho]=\beta_{hk}(\alpha)\rho\delta_{r, i} + (1-\beta_{hk}(\alpha)\rho)\delta_{h,i}, \quad i=h-1,h, h+1,\]
where $r$ and $\beta_{hk}$ are defined by:
\[r=\argmin_{j \in \{h-1, h+1\}} \{d(\theta_P, \theta_j)\},\]
\[
\beta_{hk}(\alpha)=
\begin{cases}
\alpha & \text{if} \,\,\, d(\theta_h, \theta_P) \geq \Delta\theta \\
\alpha \dfrac{d(\theta_h, \theta_P)}{\Delta\theta}& \text{if} \,\,\, d(\theta_h, \theta_P)< \Delta \theta.             
\end{cases}
 \]
We recall that $d(\cdot, \cdot)$ is defined in \eqref{eq:distance}.

%%%%%%%%%%%%%%%%%%%%%%%%
\subsection{Mathematical model}

Two last ingredients are needed before we can state the mathematical model. These are:
\begin{itemize}
\item[-] the \textit{interaction rate} with geometric features $\mu[\rho]$ : it models the frequency of interactions 
between candidate people and the walls and/or obstacles. 
If the local density is lower, it is easier for pedestrians to see the walls and doors. Thus, we set $\mu[\rho] =1-\rho$.  
\item[-] the \textit{interaction rate} with people $\eta[\rho]$: it defines the number of binary encounters per unit time. If the local density increases, then the interaction rate also increases. For simplicity, we take $\eta[\rho]= \rho$. 
\end{itemize}

The mathematical model is derived from a suitable balance of people in an elementary volume
of the space of microscopic states, considering the net flow into such volume due to
transport and interactions.
We obtain:
\begin{align}
\frac{\partial f^i}{\partial t} &+ \nabla \cdot \left( \v^i [\rho] (t, \x) f^i(t, \x) \right) \cl
& = \mathcal{J}^i[f](t, \x) \cl
& = \mathcal{J}^i_G[f](t, \x) + \mathcal{J}^i_P[f](t, \x) \cl
& = \mu[\rho] \left( \sum_{h = 1}^n \mathcal{A}_h(i) f^h(t, \x) - f^i(t,\x) \right) \cl
& \quad + \eta[\rho] \left( \sum_{h,k = 1}^n \mathcal{B}_{hk}(i) [\rho] f^h(t, \x)f^k(t, \x) - f^i(t,\x) \rho(t, \x)\right)\label{eq:model}
\end{align}
for $i= 1,2, \dots, N_d$. Functional $\mathcal{J}^i[f]$ represents the net balance of people 
that move with direction $\theta_i$ due to interactions. Since we 
consider both the interaction with the environment and with the surrounding people,
we can write $\mathcal{J}^i$ as $\mathcal{J}^i=\mathcal{J}^i_G + \mathcal{J}^i_P$, 
where $\mathcal{J}^i_G$ is an interaction between candidate people and the geometry of the environment and 
 $\mathcal{J}^i_P$ is an interaction between candidate and field people.
 
Eq.~\eqref{eq:model} is completed with eq.~\eqref{eq:rho}
for the density and eq.~\eqref{eq:v},\eqref{eq:coeff} for the velocity. 
In the next section, we will discuss a numerical method for the solution of
problem \eqref{eq:rho},\eqref{eq:v},\eqref{eq:coeff},\eqref{eq:model}.

\subsection{Full discretization} \label{sec:crowd_dyn_discr}

The approach we consider is based on a splitting method that decouples
the treatment of the transport term and the interaction term in eq.~\eqref{eq:model}.
As usual with splitting methods, the idea is to split the model into a set of subproblems that are easier to solve
and for which practical algorithms are readily available. 
Among the available operator-splitting methods, we
chose the Lie splitting scheme because it provides
a good compromise between accuracy and robustness, 
as shown in \cite{glowinski2003finite}.

Let $\Delta t>0$ be a time discretization step for the time interval $[0, T]$. 
Denote $t^k=k\Delta t$, with $k = 0, \dots, N_t$ and let $\phi^k$
be an approximation of $\phi(t^k).$ 
Given an initial condition $f^{i,0}=f^i(0, \x)$, for $i = 1, \dots, N_d$, the Lie operator-splitting scheme 
applied to problem \eqref{eq:model} reads:
For $k=0,1,2, \dots, N_t-1,$ perform the following steps:
\begin{itemize}
\item[-] {\bf Step 1}: Find $f^i$, for $i = 1, \dots, N_d$, such that\\
\begin{equation}
\begin{cases}
 \dfrac{\partial f^i}{\partial t} + \dfrac{\partial }{\partial x} \left( (v [\rho] \cos\theta_i) f^i(t, \x) \right)=0
   \,\,\, \text{on } (t^k, t^{k+1}), \\ \label{eq:step1}
f^i(t^k, \x)=f^{i,k}.
\end{cases} 
\end{equation}
 Set $f^{i,k+\frac{1}{3}}=f^i(t^{k+1}, \x)$.
 
\bigskip
\item[-] {\bf Step 2}:  Find $f^i$, for $i = 1, \dots, N_d$, such that \\
\begin{equation}
\begin{cases}
 \dfrac{\partial f^i}{\partial t} + \dfrac{\partial }{\partial y} \left( (v [\rho] \sin\theta_i) f^i(t, \x) \right)=0
   \,\,\, \text{on } (t^k, t^{k+1}),  \\ \label{eq:step2}
f^i(t^k, \x)=f^{i, k+\frac{1}{3}}. 
\end{cases}
\end{equation}
 Set $f^{i, k+\frac{2}{3}}=f^i(t^{k+1}, \x)$.

\bigskip
\item[-] {\bf Step 3}:  Find $f_i$, for $i = 1, \dots, N_d$, such that\\
\begin{equation}
\begin{cases}
 \dfrac{\partial f^i}{\partial t} = \mathcal{J}^i[f](t, \x)  \,\,\, \text{on } (t^k, t^{k+1}),  \\ \label{eq:step3}
f^i(t^k, \x)=f^{i, k+\frac{2}{3}}.
\end{cases}
\end{equation}
 Set $f^{i,k+1}=f^i(t^{k+1}, \x)$.
\end{itemize}
\bigskip

Once $f^{i, k+1}$ is computed for $i = 1, \dots, N_d$, we use eq.~\eqref{eq:rho} to get the density $\rho^{k+1}$ and
equation~\eqref{eq:v},\eqref{eq:coeff} to get the velocity magnitude at time $t^{k+1}$.

To complete the numerical method, we need to pick an appropriate numerical scheme for each subproblem.

For simplicity, we present space discretization for 
computational domain $[0, L] \times [0, H]$, with $L$ and $H$ given. 
We mesh the domain by choosing $\Delta x$ and $\Delta y$ to partition interval $[0, L]$
and $ [0, H]$, respectively. Let $N_x = L/\Delta x$ and $N_y = H/\Delta y$.
We define the discrete mesh points $\x_{pq} = (x_p, \, y_q)$ by
\begin{align}
x_p =p \Delta x~\text{with } p= 0, 1, \dots, N_x,  \quad y_q  =q \Delta y~\text{with } q= 0, 1, \dots, N_y.    \el
\end{align}
It is also useful to define 
\begin{align}
x_{p+1/2}=x_{p}+\Delta x/2=\Big(p+\frac{1}{2}\Big)\Delta x, \quad y_{q+1/2}=y_{q}+\Delta y/2=\Big(q+\frac{1}{2}\Big)\Delta y. \el
\end{align}

In order to simplify notation of the fully discrete steps 1-3, let us set $\phi = f^i$, $\theta = \theta_i$, $t_0 = t^k$, $t_f = t^{k+1}$.
Let $M$ be a positive integer ($\geq 3$, in practice). We associate with $M$ a time discretization step $\tau = (t_f - t_0)/M$
and set $t^m = t_0 + m \tau$. The fully discretized version of the Lie splitting algorithm is as follow.
\\

\noindent{\bf Discrete step 1}

\noindent Let $\phi_0 = f^{i,k}$.
Problem \eqref{eq:step1} can be rewritten as
\begin{equation}
\begin{cases}
 \dfrac{\partial \phi}{\partial t} + \dfrac{\partial }{\partial x} \left( (v [\rho] \cos\theta) \phi(t, \x) \right)=0
   \,\,\, \text{on } (t_0, t_f), \\ \label{eq:step1_bis}
\phi(t_0, \x)= \phi_0.
\end{cases} 
\end{equation}
We adopt a finite difference method that produces an approximation $\Phi_{p,q}^{m} \in \mathbb{R}$ 
of the cell average:
\[
\Phi_{p,q}^{m} \approx \dfrac{1}{\Delta x \, \Delta y} \int_{y_{q-1/2}}^{y_{q+1/2}}  \int_{x_{p-1/2}}^{x_{p+1/2}}
\phi(t^m, x, y) dx\, dy, 
\]
where $m=1, \dots, M$, $1 \leq p \leq N_x-1$ and $1 \leq q \leq N_y-1$. 
Given an initial condition $\phi_0$,  function $\phi^m$ will be approximated by $\Phi^{m}$ with
\[
\Phi^m \bigg|_{[x_{p-1/2}, \, x_{p+1/2}] \times [y_{q-1/2}, \, y_{q+1/2}]} = \Phi_{p,q}^{m}
\]

The Lax-Friedrichs method for problem \eqref{eq:step1_bis} can be written in conservative form as follows:
\[\Phi_{p,q}^{m+1}=\Phi_{p,q}^{m}- \dfrac{\tau}{\Delta x}\Big(  \mathcal{F}(\Phi_{p,q}^{m}, \Phi_{p+1,q}^{m}) 
- \mathcal{F}(\Phi_{p-1,q}^{m}, \Phi_{p,q}^{m}) \Big)\]
where
\[ \mathcal{F}(\Phi_{p,q}^{m}, \Phi_{p+1,q}^{m}) =\dfrac{\Delta x}{2\tau}(\Phi_{p,q}^{m}- \Phi_{p+1,q}^{m}) + \dfrac{1}{2} \Big( (v [\rho^m_{p,q}] \cos\theta) \Phi_{p,q}^{m}+(v [\rho_{p+1,q}^m] \cos\theta) \Phi_{p+1,q}^{m} \Big).  \]
\bigskip

\noindent{\bf Discrete step 2}

\noindent Let  $\phi_0 = f^{i,k+\frac{1}{3}}$.
Problem \eqref{eq:step2} can be rewritten as
\begin{equation}
\begin{cases}
 \dfrac{\partial \phi}{\partial t} + \dfrac{\partial }{\partial y} \left( (v [\rho] \sin\theta) \phi(t, \x) \right)=0
   \,\,\, \text{on } (t_0, t_f), \cl
\phi(t_0, \x)= \phi_0 \cl
\end{cases}
\end{equation}
Similarly to step 1, we use the conservative Lax-Friedrichs scheme:
\[\Phi_{p,q}^{m+1}=\Phi_{p,q}^{m}- \dfrac{\tau}{\Delta y}\Big(  \mathcal{F}(\Phi_{p,q}^{m}, \Phi_{p,q+1}^{m}) 
- \mathcal{F}(\Phi_{p,q-1}^{m} \Phi_{p,q}^{m}) \Big)\]
where
\[ \mathcal{F}(\Phi_{p,q}^{m}, \Phi_{p,q+1}^{m}) =\dfrac{\Delta y}{2\tau}(\Phi_{p,q}^{m}-\Phi_{p,q+1}^{m}) + \dfrac{1}{2} \Big( (v [\rho^m_{p,q}] \sin\theta) \Phi_{p,q}^{m}+(v [\rho_{p,q+1}^m] \sin\theta) \Phi_{p,q+1}^{m} \Big).  \]
\bigskip

\noindent{\bf Discrete step 3}

\noindent Let $\mathcal{J} = \mathcal{J}^i $ and $\phi_0 = f^{i,k+\frac{2}{3}}$.
Problem \eqref{eq:step3} can be rewritten as
\begin{equation}
\begin{cases}
 \dfrac{\partial \phi}{\partial t} = \mathcal{J}[f](t, \x)  \,\,\, \text{on } (t_0, t_{f}), \cl
\phi(t_0, \x)= \phi_0. \el
\end{cases}
\end{equation}
For the approximation of the above problem, we use the Forward Euler scheme:
\[ \Phi_{p,q}^{m+1}= \Phi_{p,q}^{m} + \tau \Big (\mathcal{J}^m[F^{m}] \Big ), \] 
where $F^m$ is the approximation of the reduced distribution function \eqref{eq:f} at time $t^m$.
\medskip

For stability, the subtime step $\tau$ is chosen to satisfy the Courant-Friedrichs-Lewy (CFL) 
condition (see, e.g., \cite{leveque1992numerical}):
\[\max \Bigg\{\cfrac{\tau}{\Delta x}, \,\, \cfrac{\tau}{\Delta y} \Bigg\}\leq 1. \]

\subsection{Numerical results} \label{sec:crowd_dyn_num_res}

\begin{figure}[htbp]
\centering
\begin{overpic}[width=0.39\textwidth, grid=false]{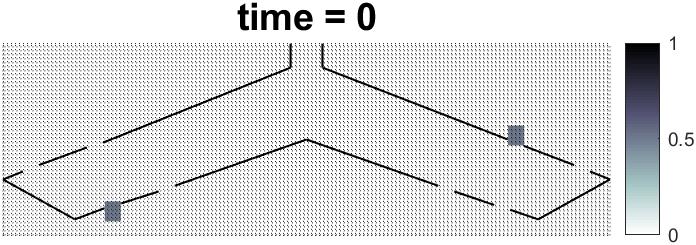}
\put(16.2, 4.8){\vector(0,1){5}}
\put(74.2, 15.9){\vector(0,-1){5}}
\end{overpic}
\begin{overpic}[width=0.39\textwidth, grid=false]{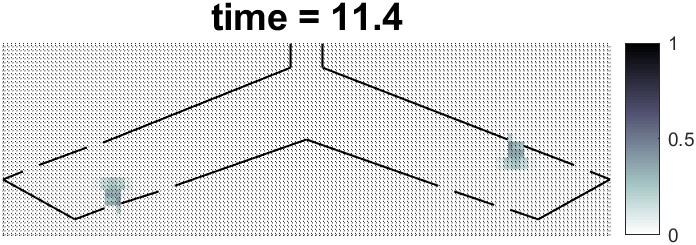}
\end{overpic}
\begin{overpic}[width=0.39\textwidth, grid=false]{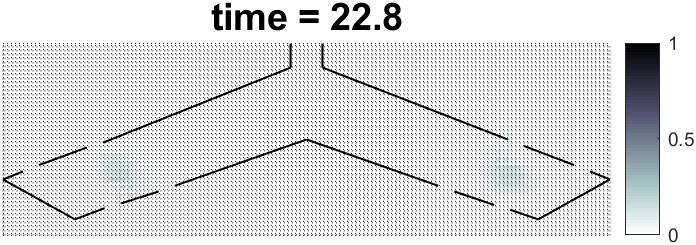}
\end{overpic}
%\begin{overpic}[width=0.49\textwidth, grid=false]{./test1_airport_time_6.jpg}
%\end{overpic}
\begin{overpic}[width=0.39\textwidth, grid=false]{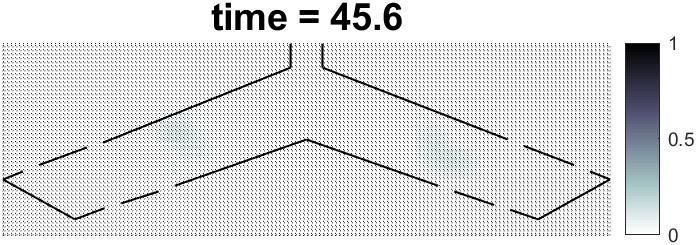}
\end{overpic}
\begin{overpic}[width=0.39\textwidth, grid=false]{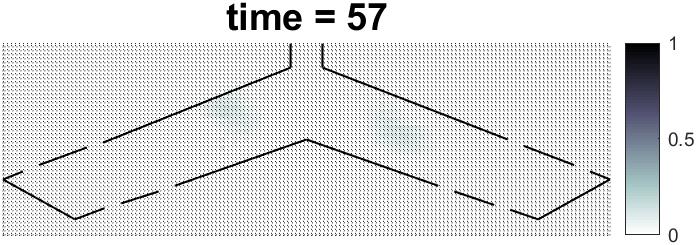}
\end{overpic}
%\begin{overpic}[width=0.49\textwidth, grid=false]{./test1_airport_time_12.jpg}
%\end{overpic}
\begin{overpic}[width=0.39\textwidth, grid=false]{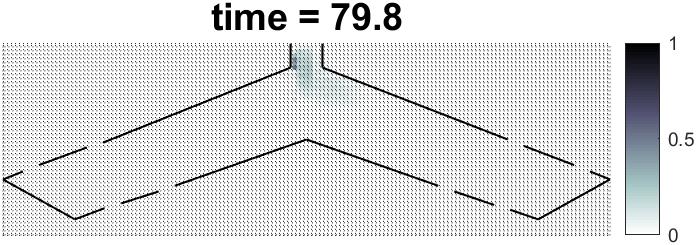}
\end{overpic}
\caption{Test 1$a$: Evacuation process of 404 people grouped into two clusters with initial directions
$\theta_3$ (group in the left wing) and $\theta_7$ (group in the right wing).}
\label{airport}
\end{figure}

\begin{figure}[htb]
\centering
\begin{overpic}[width=0.39\textwidth, grid=false]{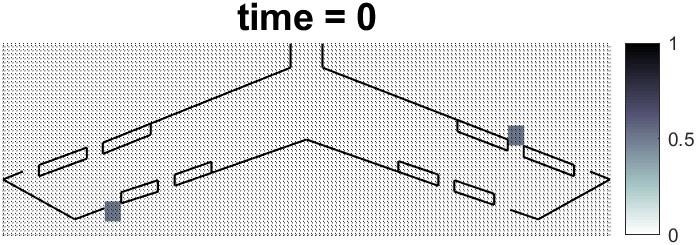}
\put(16.2, 4.8){\vector(0,1){5}}
\put(74.2, 15.9){\vector(0,-1){5}}
\end{overpic}
\begin{overpic}[width=0.39\textwidth, grid=false]{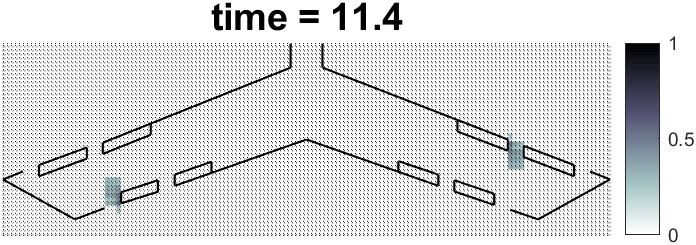}
\end{overpic}
\begin{overpic}[width=0.39\textwidth, grid=false]{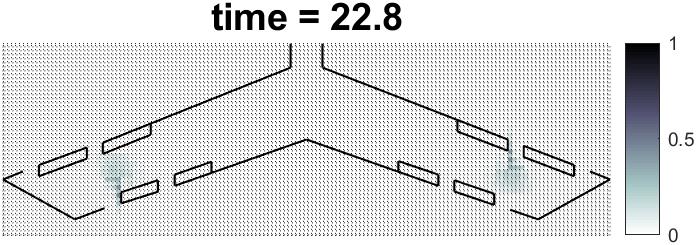}
\end{overpic}
%\begin{overpic}[width=0.49\textwidth, grid=false]{./test2_airport_time_6.jpg}
%\end{overpic}
\begin{overpic}[width=0.39\textwidth, grid=false]{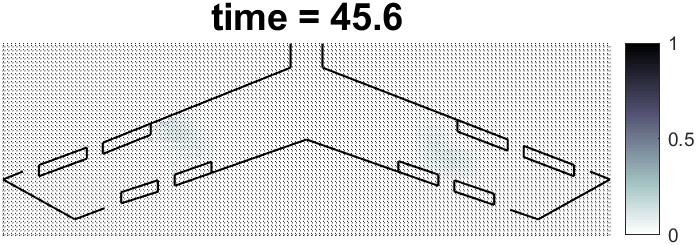}
\end{overpic}
\begin{overpic}[width=0.39\textwidth, grid=false]{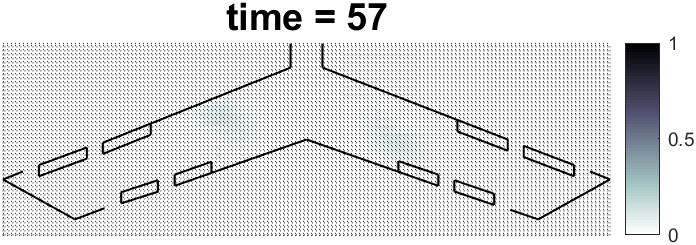}
\end{overpic}
%\begin{overpic}[width=0.49\textwidth, grid=false]{./test2_airport_time_12.jpg}
%\end{overpic}
\begin{overpic}[width=0.39\textwidth, grid=false]{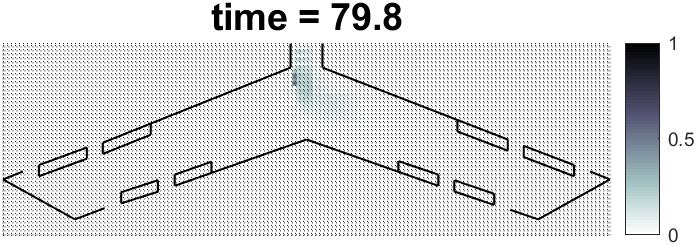}
\end{overpic}
\caption{Test 1$b$: Evacuation process of 404 people grouped into two clusters with initial directions
$\theta_3$ (group in the left wing) and $\theta_7$ (group in the right wing).}
\label{airport_chair}
\end{figure}

\begin{figure}[htb]
\centering
\begin{overpic}[width=0.39\textwidth, grid=false]{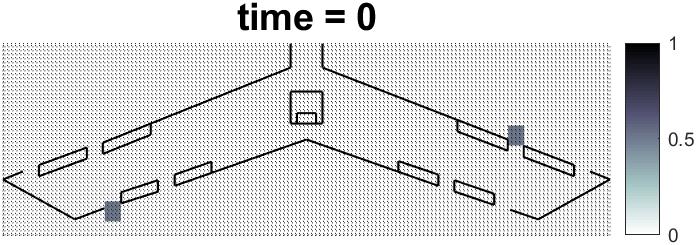}
\put(16.2, 4.8){\vector(0,1){5}}
\put(74.2, 15.9){\vector(0,-1){5}}
\end{overpic}
\begin{overpic}[width=0.39\textwidth, grid=false]{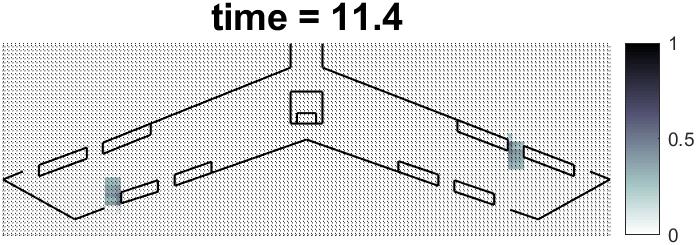}
\end{overpic}
\begin{overpic}[width=0.39\textwidth, grid=false]{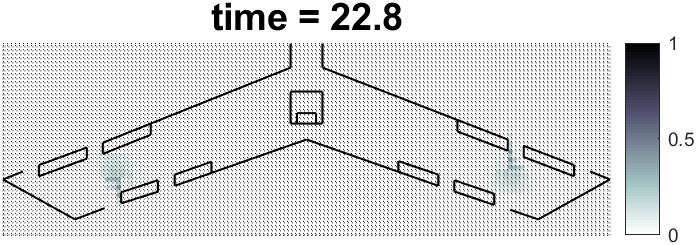}
\end{overpic}
%\begin{overpic}[width=0.49\textwidth, grid=false]{./test3_airport_time_6.jpg}
%\end{overpic}
\begin{overpic}[width=0.39\textwidth, grid=false]{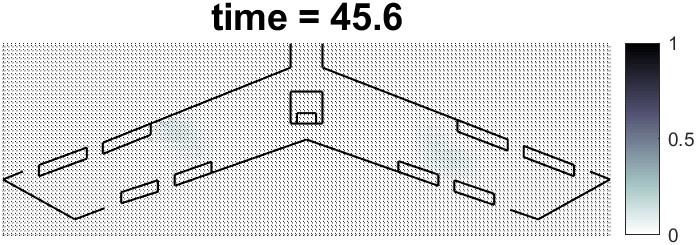}
\end{overpic}
\begin{overpic}[width=0.39\textwidth, grid=false]{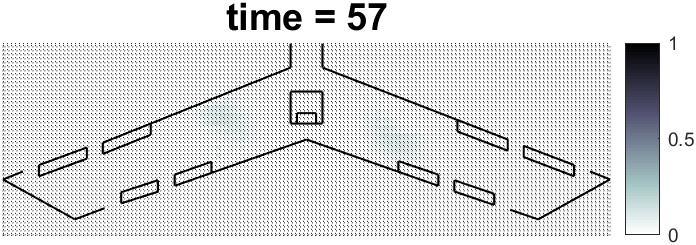}
\end{overpic}
%\begin{overpic}[width=0.49\textwidth, grid=false]{./test3_airport_time_12.jpg}
%\end{overpic}
\begin{overpic}[width=0.39\textwidth, grid=false]{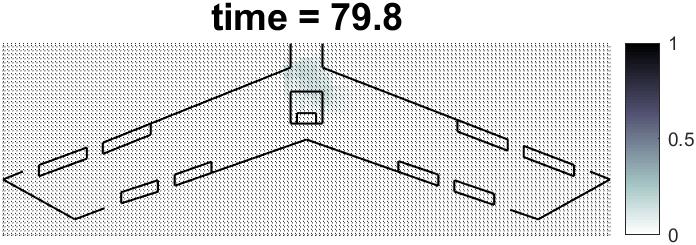}
\end{overpic}
\caption{Test 1$c$: Evacuation process of 404 people grouped into two clusters with initial directions
$\theta_3$ (group in the left wing) and $\theta_7$ (group in the right wing).}
\label{airport_piano}
\end{figure}

We consider a part of Houston's William P. Hobby Airport as the walkable domain.
The terminal has an upside down V shape with eight 
gates (four per wing) and an entrance/exit at the top, which is 
$6.8$ m wide. See Fig.~\ref{airport} (top left panel). 
The shape and the size of the terminal (each wing is about 136 m long and 20 m wide)
are realistic, while the number of gates is reduced 
for simplicity. We consider the following geometries:
\begin{itemize}
\item[-] Configuration $a$: no obstacle in the terminal. See Fig.~\ref{airport} (top left panel). 
\item[-] Configuration $b$: waiting area chairs are located near each terminal. See Fig.~\ref{airport_chair} (top left panel). 
\item[-] Configuration $c$: in addition to the waiting area chairs, a large obstacle, like a temporary store, is located at the intersection of the two wings. See Fig.~\ref{airport_piano} (top left panel). 
\end{itemize}

In these configurations, we run two sets of simulations: 
\begin{itemize}
\item[-] Test 1: a total of 404 passengers from two planes at the opposite 
ends of the terminal walk through the terminal to reach the exit.
\item[-] Test 2: the 404 passengers have the same target as in test 1 but there is
an additional group of 202 passengers that enter the terminal through the entrance and are directed to a gate.
\end{itemize}
The aim is to compute the egress time, i.e. the total time it takes all the passengers the leave the
terminal through either the exit or a gate. 

For all the simulations, we consider eight different velocity directions $N_d=8$ in the discrete set:
\[ I_{\theta}=\left \{ \theta_{i}= \frac{i-1}{8} 2\pi : i = 1, \dots, 8 \right \}. \]
In order to work with dimensionless quantities, we define the following reference quantities:
$D = 137.5$ m, $V_M = 2$ m/s,  and $\rho_M = 7$ people/m$^2$.
Once the results are computed, we convert them back to dimensional quantities.

%The computational domain in the xy-plane is approximately $[0, 288.8] \times [0, 91.2]$ with 4 exits located at the end of each wing. 

\begin{figure}[htb]
\centering
\begin{overpic}[width=0.45\textwidth, grid=false]{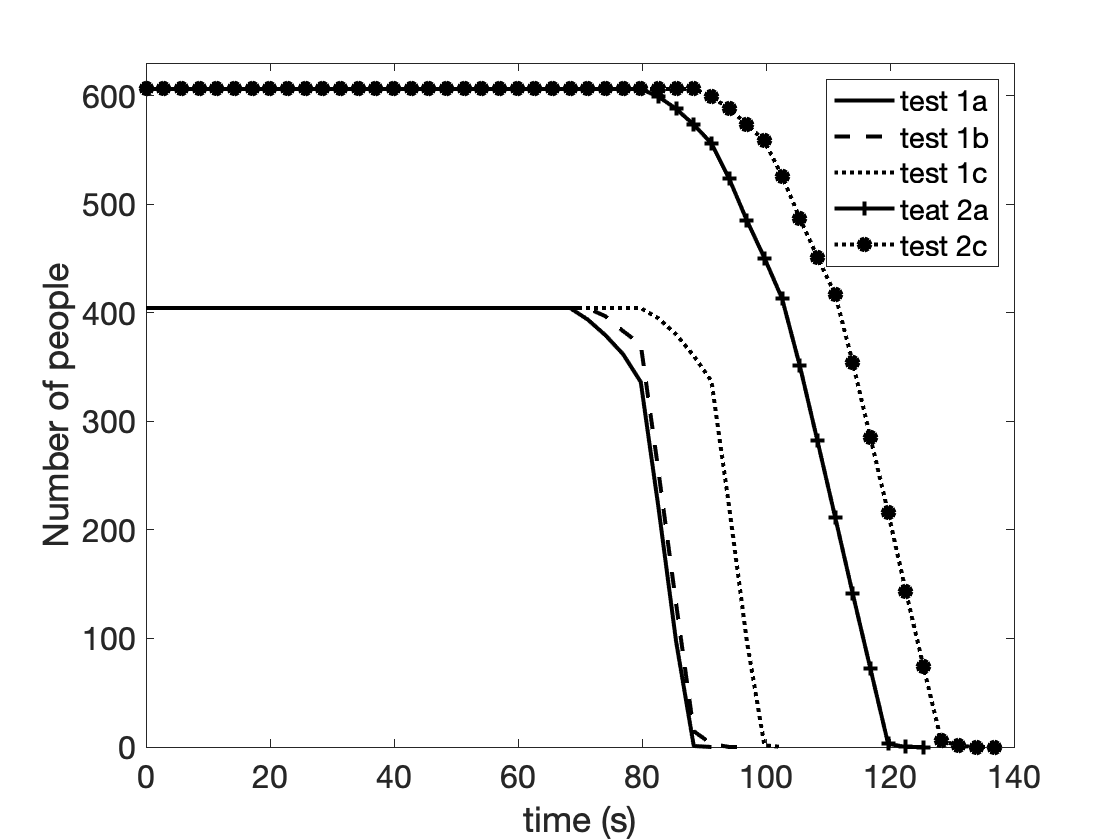}
\end{overpic}
\caption{Number of passengers inside the terminal over time for tests 1$a$, 1$b$, 1$c$, 
2$a$, and 2$c$.} \label{fig:people_terminal}
\end{figure}

We consider a mesh with $\Delta x = \Delta y = 1.9$ m. The time step is set to $\Delta t = 5.7$ s
and we choose $M = 3$. 
Fig. \ref{airport}, \ref{airport_chair}, and \ref{airport_piano} show the density
computed at different times for tests 1$a$, 1$b$, and 1$c$, respectively.
For the large obstacle in configuration $c$,  we use an effective area that is a square with side 15.2 m,
 while the actual obstacle is a rectangle with dimensions of 9.5 m in length and 4.75 m in width.
 The reader interested in learning more about how obstacles are handled is referred to \cite{kim_quaini}.
 In configuration $a$, we observe a denser crowd only when the several passengers 
 reach the exit, as shown in Fig.~\ref{airport} (lower right panel).
 Configuration $b$ creates dense gatherings also when passengers deplane and
 their motion is restricted by the waiting area chairs. See Fig.~\ref{airport_chair}
 for times $t = 11.4, 22.8$ s.
Nonetheless, we observe a similar evacuation dynamics between configurations
$a$ and $b$, indicating that the waiting area chairs do not hinder the evacuation 
process. Compare Fig.~\ref{airport} with Fig.~\ref{airport_chair}. 
This is confirmed by Fig.~\ref{fig:people_terminal}, which shows
the number of passengers inside the terminal over time for all the tests. 
The curves for tests $1a$ and $1b$ are either superimposed or very close
to each other over the entire time interval. On the other hand, we see
that the presence of a large obstacle at the intersection of the two terminal wings
increases the egress time by over 10 s. Also compare Fig.~\ref{airport_chair} with Fig.~\ref{airport_piano}, 
bottom right panels.

\begin{figure}[htb]
\centering
\begin{overpic}[width=0.39\textwidth, grid=false]{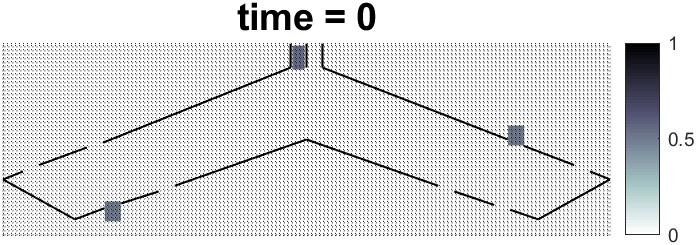}
\put(16.2, 4.8){\vector(0,1){5}}
\put(74.2, 15.9){\vector(0,-1){5}}
\put(43, 27){\vector(0,-1){5}}
\end{overpic}
\begin{overpic}[width=0.39\textwidth, grid=false]{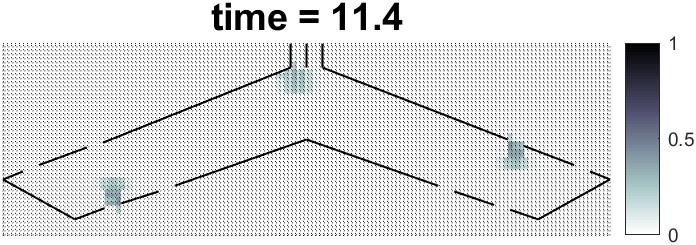}
\end{overpic}
\begin{overpic}[width=0.39\textwidth, grid=false]{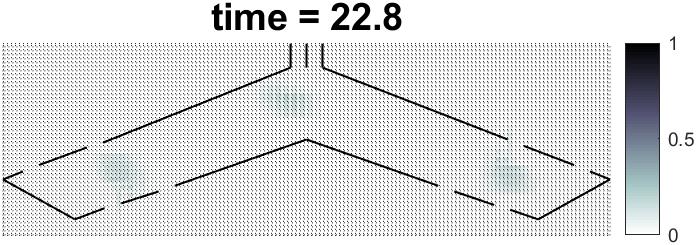}
\end{overpic}
%\begin{overpic}[width=0.49\textwidth, grid=false]{./test4_airport_time_6.jpg}
%\end{overpic}
\begin{overpic}[width=0.39\textwidth, grid=false]{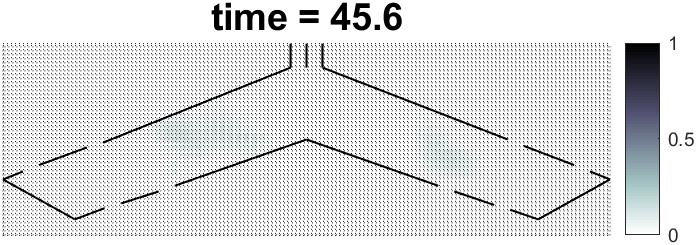}
\end{overpic}
\begin{overpic}[width=0.39\textwidth, grid=false]{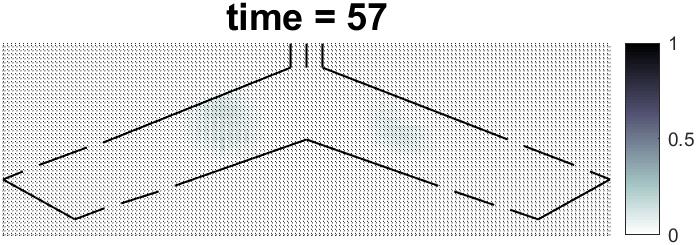}
\end{overpic}
%\begin{overpic}[width=0.49\textwidth, grid=false]{./test4_airport_time_12.jpg}
%\end{overpic}
\begin{overpic}[width=0.39\textwidth, grid=false]{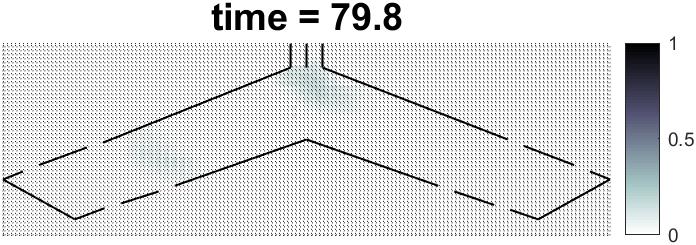}
\end{overpic}
\begin{overpic}[width=0.39\textwidth, grid=false]{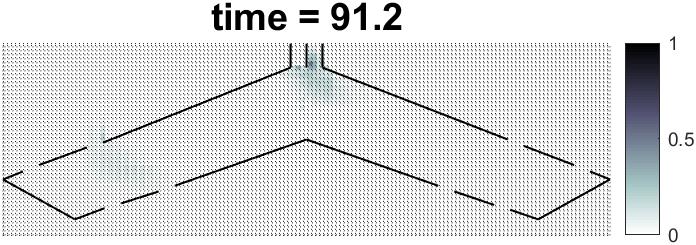}
\end{overpic}
\begin{overpic}[width=0.39\textwidth, grid=false]{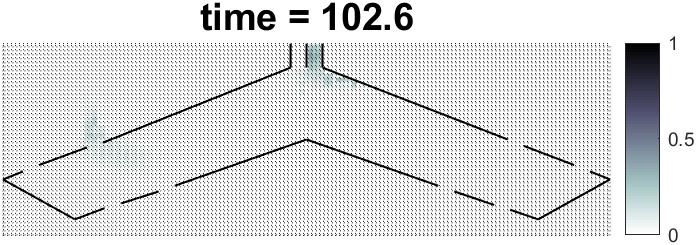}
\end{overpic}
\caption{Test 2$a$: Evacuation process of 404 people grouped into two clusters with initial directions
$\theta_3$ (group in the left wing) and $\theta_7$ (group in the right wing) at the same time as
a third group of 202 people with initial direction $\theta_7$
enters the airport and is directed to a gate in the left wing. }
\label{airport_two}
\end{figure}

\begin{figure}[htb]
\centering
\begin{overpic}[width=0.39\textwidth, grid=false]{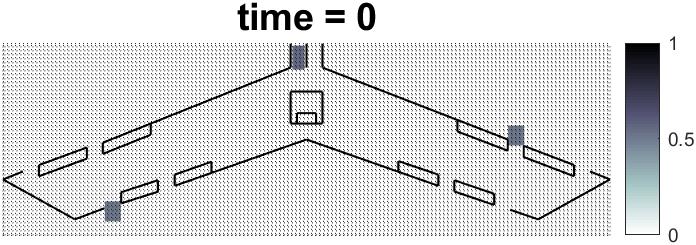}
\put(16.2, 4.8){\vector(0,1){5}}
\put(74.2, 15.9){\vector(0,-1){5}}
\put(43, 27){\vector(0,-1){5}}
\end{overpic}
\begin{overpic}[width=0.39\textwidth, grid=false]{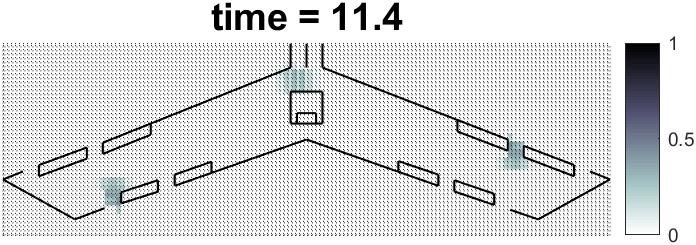}
\end{overpic}
\begin{overpic}[width=0.39\textwidth, grid=false]{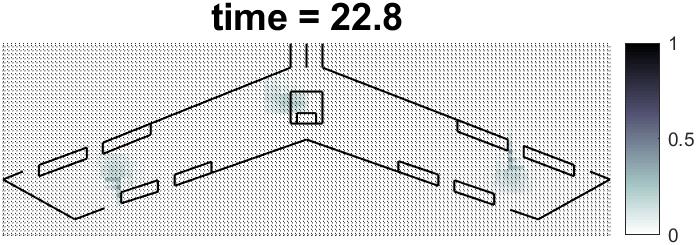}
\end{overpic}
%\begin{overpic}[width=0.49\textwidth, grid=false]{./test4_airport_time_6.jpg}
%\end{overpic}
\begin{overpic}[width=0.39\textwidth, grid=false]{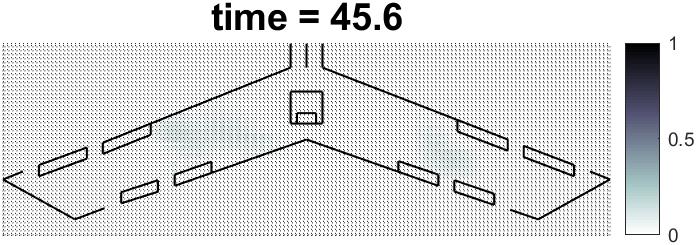}
\end{overpic}
\begin{overpic}[width=0.39\textwidth, grid=false]{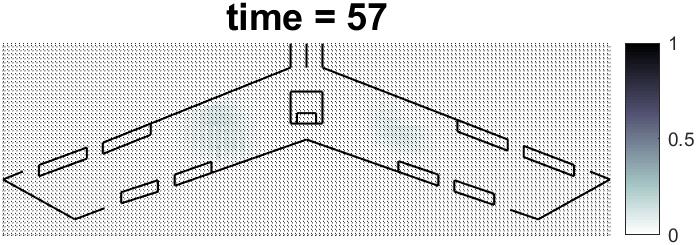}
\end{overpic}
%\begin{overpic}[width=0.49\textwidth, grid=false]{./test5_airport_time_12.jpg}
%\end{overpic}
\begin{overpic}[width=0.39\textwidth, grid=false]{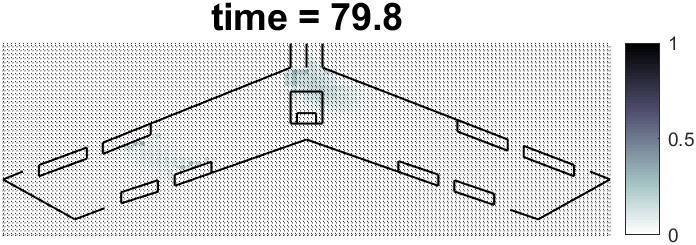}
\end{overpic}
\begin{overpic}[width=0.39\textwidth, grid=false]{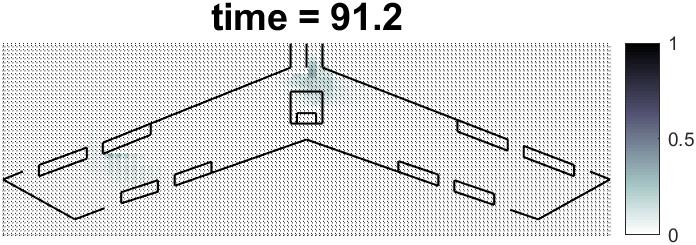}
\end{overpic}
\begin{overpic}[width=0.39\textwidth, grid=false]{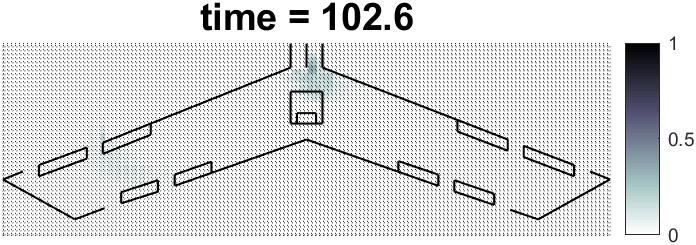}
\end{overpic}
\caption{Test 2$c$: Evacuation process of 404 people grouped into two clusters with initial directions
$\theta_3$ (group in the left wing) and $\theta_7$ (group in the right wing) at the same time as
a third group of 202 people with initial direction $\theta_7$
enters the airport and is directed to a gate in the left wing.}
\label{airport_piano_two}
\end{figure}

Because of the similarities in the evacuation process for tests 1$a$ and 1$b$, 
we decided to run test 2 only in configurations $a$ and $c$. The density
computed at different times for these two tests are shown in Fig.~\ref{airport_two} and \ref{airport_piano_two}.
In test 2 the exit size is halved because half of the top corridor is used as an entrance. From Fig.~\ref{airport_two} and \ref{airport_piano_two}
(second row, right panel), we see that by time $t = 45.6$ s the two groups of
passengers with opposite directions (heading to the exit vs to the gate) have met. As expected, halving 
the exit size leads to a longer evacuation process (for example, compare Fig.~\ref{airport}
and \ref{airport_two}) and creates a dense crowd at the exit. 
See Fig.~\ref{airport_two} and \ref{airport_piano_two} (bottom right panel).
From Fig.~\ref{fig:people_terminal}, we see that the increase in ingress
time from test 1$a$ to 2$a$ is about 30 s, while it is about 40 s from test 1$c$ to 2$c$.

The results presented in this subsection corroborate the effectiveness of some strategies 
adopted in airports during the COVID-19 pandemic: dedicated, distant sites for
entrances and exits, and minimization of the obstacles inside the terminal. These strategies
are conducive to short egress times and limit congregation points, thereby containing the spreading
of COVID-19.

\section{Contagion model in one dimension}\label{sec:contagion}
%%%%%%%%%%%%%%%%%%%%%%%%%%%%

We start from an agent-based model at the microscopic level. 
We consider a group of $N$ people, $N_h$ of whom are healthy or not spreading the disease yet, 
while the remaining $N_s = N - N_h$ are in the spreading phase of the disease.
If person $n$ belongs to the former group, we denote with $q_n \in [0, 1)$
their level of exposure to people spreading the disease, with the underlying
idea that the more a person is exposed the more likely they are to get infected.
If person $n$ belongs to the latter group, then $q_n = 1$ and this value stays constant throughout the entire simulation time. 
In addition, let $x_n(t)$ and $v_n(t)$ denote the position and speed of person $n$.

The microscopic model reads for $n= 1,2,3, \dots, N$:
\begin{align}
&\frac{dx_n}{dt}=v_n\cos\theta_n, ~ \frac{dq_n}{dt}=\gamma \max\{(q_n^*-q_n), 0\},~ q_n^* =\frac{\sum_{m=1}^{N} \kappa_{n,m}q_m}{\sum_{m=1}^{N}\kappa_{n,m}}, \label{eq:diseasef}
\end{align}
where the walking speed $v_n$ and walking direction $\theta_n$ are given.
In the future, we will combine the model in this section with the model presented in 
Sec.~\ref{sec:ped_model} that will provide walking speed and direction. 
In model \eqref{eq:diseasef}, $q_n^{\ast}$
corresponds to a weighted average ``level of sickness'' surrounding person $n$, with 
$\kappa_{n,m}$ that serves as the weight in the average. 
We define $\kappa_{n,m}$ as follow
\begin{equation}\label{kappa}
\kappa_{n,m} = \kappa(|x_n-x_m|) = \dfrac{R}{(|x_n-x_m|^2+R^2) \pi} 
\end{equation}
Notice that the interaction kernel 
is a decreasing function of mutual distance between two people and is parametrized by an interaction distance $R$,
set so that the value of $\kappa_{n,m}$ is ``small'' at about 6 ft or 2 m.
Parameter $\gamma$ in \eqref{eq:diseasef} describes the contagion interaction strength: 
for $\gamma = 0$ there is no contagion, while for
$\gamma \neq 0$ the contagion is faster the larger the value of $\gamma$. 
Note that obviously the
level of exposure can only increase over time. The second equation in \eqref{eq:diseasef}
also ensure that the people spreading the disease will constantly have $q_n = 1$ in time.

From the agent-based model \eqref{eq:diseasef}, we derive a model at the kinetic level. 
Denote the empirical distribution by
\begin{equation}
h^N=\frac{1}{N}\sum_{n=1}^N\delta(x-x_n(t))\delta(q-q_n(t)), \cl
\end{equation}
where $\delta$ is the Dirac delta measure.
We assume that the people remain in a fixed compact domain $(x_n(t), q_n(t)) \in \Omega \subset \mathbb{R}^2$ for all $n$ and for the entire time interval under consideration.
Prohorov's theorem implies that the sequence $\{h^N\}$ is relatively compact in the weak$^{\ast}$ sense. 
Therefore, there exists a subsequence $\{h^{N_k}\}_k$ such that $h^{N_k}$ converges to $h$ with weak$^{\ast}$-convergence in $\mathcal{P}(\mathbb{R}^2$) and pointwise convergence in time as $k \rightarrow \infty$. 
Here, $\mathcal{P}(\mathbb{R}^2$) denotes the space of probability measures on $\mathbb{R}^2$.
%\WOcomment{What is the regularity of the probability measures in this space? Is this the space of Borel probability measures, or something else?} 

Let $\psi \in C_0^1(\mathbb{R}^2)$ be a test function. We have
\begin{align}\label{eq:testfunction}
\frac{d}{dt} \langle h^N, \psi \rangle_{x,q} 
&=\frac{d}{dt} \biggl< \frac{1}{N}\sum_{n=1}^{N}\delta(x-x_n(t))\delta(q-q_n(t)), \psi \biggr>_{x,q} \cl
&=\frac{d}{dt} \frac{1}{N}\sum_{n=1}^{N}\psi(x_n(t), q_n(t)) \cl
&=\frac{1}{N}\sum_{n=1}^{N} \left( \psi_x v_n \cos \theta_n + \psi_q\gamma\max\{(q_n^*-q_n), 0\} \right) \cl
&= \langle h^N, \psi_x v \cos \theta_n \rangle_{x,q} + \frac{\gamma}{N}\sum_{n=1}^{N}\psi_q \max\left\{\left( \frac{\sum_{m=1}^{N} \kappa_{n,m}q_n}{\sum_{m=1}^{N}\kappa_{n,m}}- q_n \right), 0 \right\},
\end{align} 
where $\langle \cdot \rangle_{x,q}$ means integration against both $x$ and $q$. 

Let us define
\begin{equation}
\rho(x)=\frac{1}{N}\sum_{n=1}^{N}\delta(x-x_n) \cl
\end{equation}
and
\begin{equation}
m(x)=\biggl< q, \frac{1}{N}\sum_{m=1}^{N}\delta(x-x_m)\delta(q-q_m) \biggr>_{x,q}
=\frac{1}{N}\sum_{m=1}^{N}\delta(x-x_m)q_m, \cl
\end{equation}
We have
\begin{align}
\frac{1}{N}\sum_{m=1}^{N}\kappa(|x_n-x_m|)&= \biggl< \kappa(|x_n-\tilde{x}|), \frac{1}{N}\sum_{m=1}^{N}\delta(\tilde{x}-x_m) \biggr>_{x}
=\kappa \ast\rho(x_n), \cl
\frac{1}{N}\sum_{m=1}^{N}\kappa(|x_n-x_m|)q_m&= \biggl< \kappa(|x_n-\tilde{x}|), \frac{1}{N}\sum_{m=1}^{N}\delta(\tilde{x}-x_m)q_m \biggr>_{x}
=\kappa \ast m(x_n), \el
\end{align}
where $\langle \cdot \rangle_x$ means integration only in $x$.
Then, we can rewrite eq.~($\ref{eq:testfunction}$) as
\begin{equation}\label{eq3_rewritten}
\frac{d}{dt} \langle h^{N}, \psi \rangle_{x,q}= \langle h^N, \psi_x v \cos \theta \rangle_{x,q} + \gamma \biggl< h^{N}, \psi_{q} \max \left\{ \frac{\kappa \ast m}{\kappa \ast \rho}-q, 0 \right\} \biggr>_{x,q}.
\end{equation}
Via integration by parts, eq.~\eqref{eq3_rewritten} leads to
\begin{equation}\label{eq:Nkineticsystem}
h_t^N+ (v\cos\theta~ h^N)_x+\gamma(\max\{(q^\ast-q),0\}h^N)_q = 0, 
\end{equation}
where $q^*$ is the local \emph{average} sickness level weighted by \eqref{kappa}:
\begin{equation} \label{q_act}
q^{\ast}(t,x)= \frac{\iint\kappa(|x-\overline{x}|) h(t,\overline{x},q)qdqd\overline{x}}{\iint \kappa(|x-\overline{x}|)h(t,\overline{x},q)dqd\overline{x}}.
\end{equation}
Sick people that are in the spreading phase of the disease weight more in the average since they have the highest value of
$q$, nonetheless exposed people contribute to the average level of sickness too
since they might spread the virus they recently got exposed to (recall we
are simulating short periods of time), e.g., by close contact. 

Now letting $k \rightarrow \infty$, the subsequence $h^{N_k}$ formally leads to the limiting kinetic equation
\begin{equation} \label{eq:2kineticeq}
h_t+ (v \cos\theta~h)_x+\gamma(\max \{(q^\ast-q), 0\} h)_q = 0, 
\end{equation}
where $h(t,x,q)$ is the probability of finding at time $t$ and position $\x$ 
a person with level of exposure $q$ if $q \in [0, 1)$
or a person spreading the disease if $q=1$.  

Finally, we note that while modeling motion and disease spreading in one dimension (spatial variable $x$), 
eq.~\eqref{eq:2kineticeq} is a 2D problem in variables $x$ and $q$. 
Modeling pedestrian motion in two dimensions would lead to 
a 3D problem that requires a carefully designed numerical scheme to contain the computational costs. 
This is currently under investigation.

%%%%%%
\subsection{Full discretization}\label{sec:contagion_full_disc}
%%%%%%
We present a space and time discretization for eq.~(\ref{eq:2kineticeq}). Let $x \in [0, D]$ and $q \in [0, 1]$.
Given $N_x = D/\Delta x$, the discrete mesh points $x_{p}$ are given by
\begin{align}
x_p =p \Delta x, \quad x_{p+1/2}=x_{p}+ \frac{\Delta x}{2}=\Big(p+\frac{1}{2}\Big)\Delta x, \label{eq:x_p}
\end{align}
for $p= 0, 1, \dots, N_x$.
We partition $[0, 1]$ into subintervals $[q_{l-\frac{1}{2}}, q_{l+\frac{1}{2}}]$, with $l \in 1,2,\dots, N_q$, where 
\begin{align}
q_l =l \Delta q, \quad q_{l+1/2}=q_{l}+ \frac{\Delta q}{2}=\Big(l+\frac{1}{2}\Big)\Delta q. \el
\end{align}
For simplicity, we assume that all subintervals have equal length $\Delta q$.
The two partitions induce a partition of domain $[0, D] \times [0, 1]$ into cells.
The time step $\Delta t$ is chosen as 
\[ \Delta t \leq \min\Bigg\{ \frac{\Delta x}{\max_{p} v_p}, \frac{\Delta q}{2\gamma \max_{l}q_{l}} \Bigg\}\] 
to satisfy the Courant-Friedrichs-Lewy (CFL) condition.

Let us denote $h_{j, l}=h(t, x_j, q_l)$ and $q_j^ \ast=q^\ast(t, x_j)$.
We consider a first-order semi-discrete upwind scheme for eq.~\eqref{eq:2kineticeq} 
adapted from one of the methods used in Ref.~\cite{Bertozzi2015}, which reads:
\begin{equation} \label{eq:disease_modified}
\partial_t h_{j, l}+\frac{\eta_{j, l}-\eta_{j-1, l}}{\Delta x}
+\gamma \frac{\xi_{j, l+\frac{1}{2}}-\xi_{j, l-\frac{1}{2}}}{\Delta q} = 0,
\end{equation}
where
\begin{align}
\eta_{j, l} & = v_j \cos \theta_j~h_{j, l},\cl
\xi_{j, l+\frac{1}{2}} & = \max\left\{\left(q^\ast_j-q_{l+\frac{1}{2}}\right), 0\right\}  h_{j, l}. \el
\end{align}
For the time discretization of problem~\eqref{eq:disease_modified}, we use the forward Euler scheme:
\begin{equation}
 h_{j,l}^{m+1} =  h_{j,l}^{m}-\Delta t \Bigg ( \cfrac{\eta^{m}_{j, l}-\eta^{m}_{j-1, l}}{\Delta x}
+\gamma \cfrac{\xi^{m}_{j, l+\frac{1}{2}}-\xi^{m}_{j, l-\frac{1}{2}}}{\Delta q} \Bigg ).
\end{equation}

The discretization scheme in this section is only first order in space and time. The numerical 
errors are expected to introduce significant dissipation in the numerical solution. Extension to higher order
discretization schemes is possible (see, e.g., \cite{kim_quaini2020b,kim_quaini2020,Bertozzi2015})
but will not be considered for this paper.

%%%%%%%
\subsection{Numerical results}\label{sec:num_res_contagion}
 %%%%%
 
\begin{figure}[htb]
\centering
\begin{overpic}[width=0.45\textwidth, grid=false]{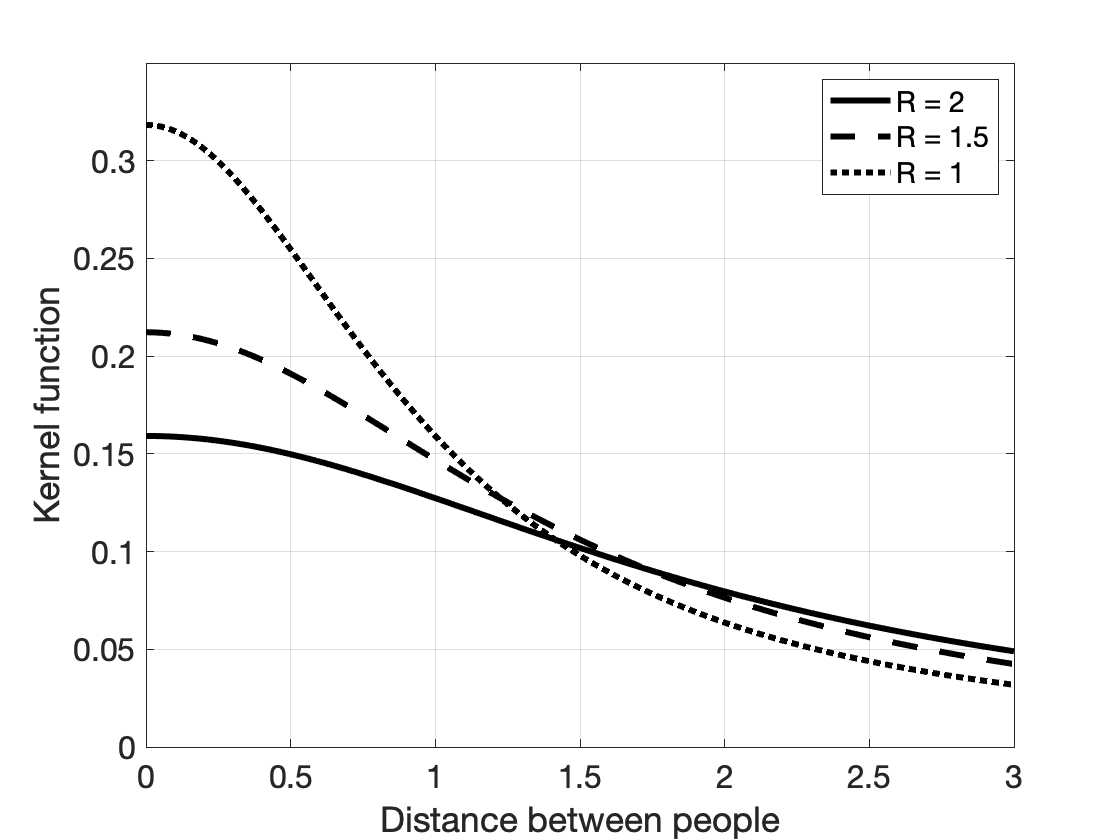}
\end{overpic}
\caption{Kernel function vs the distance between people for interaction radius $R = 1, 1.5, 2$.} \label{fig:kappa}
\end{figure}

We test the approach presented in Sec.~\ref{sec:contagion_full_disc} on a series of 1D problems, corresponding
to unidirectional pedestrian flow in a narrow corridor. For all the problems, 
the computational domain in the $xq$-plane is $[0, 10] \times [0, 1]$ and it is occupied by a group of 
40 people. We set $R = 1$ m since this
choice makes the value of the kernel function relatively small at a distance of 2 m (or about 6 ft). See Fig.~\ref{fig:kappa}.
The dimensionless quantities are obtained 
by using the following reference quantities: $D=10$ m, $V_M= 1$ m/s, $T = 10$ s, $\rho_M = 4$ people/m.
In all the tests, we take the initial density to be constant in space
and equal to $\rho_M$.

We take $\Delta x = 0.1$ m and $\Delta q = 0.01$. We will consider two values for contagion strength $\gamma = 100$
and $\gamma = 50$, with the associated respective time steps $\Delta t=0.00005$ s and  $\Delta t=0.0001$ s.
First, we keep the group of people still (i.e., $v = 0$) to observe how the level of exposure to the disease evolve.
Then, in a second set of tests, we change to $v = 1$ m/s and see how the motion affects the spreading. 
We run each simulation for $t \in (0, 10]$ s.

{\bf Tests with $\boldsymbol{v = 0}$.} We consider two initial conditions
\begin{itemize}
\item[-] IC1: people that are certainly spreading (i.e., $q = 1$) are located at $x \in [0, 4]$ m and $x \in [6,10]$ m, while
in $x \in (4,6)$ m we place people that have certainly not been exposed (i.e, $q = 0$). 
\item[-] IC2: people that are certainly spreading (i.e., $q = 1$) are located at $x \in [0, 2]$ m and $x \in [8,10]$ m, while
the rest of the people located  in $x \in (2,8)$ m have certainly not been exposed (i.e, $q = 0$). 
\end{itemize}
All the healthy people in IC1 are exposed to both groups of spreading people, while in IC2 some 
healthy people are exposed to one group of spreading people and the centrally located
healthy people are not exposed.

\begin{figure}[htb]
\centering
\begin{overpic}[width=0.25\textwidth,grid=false,tics=10]{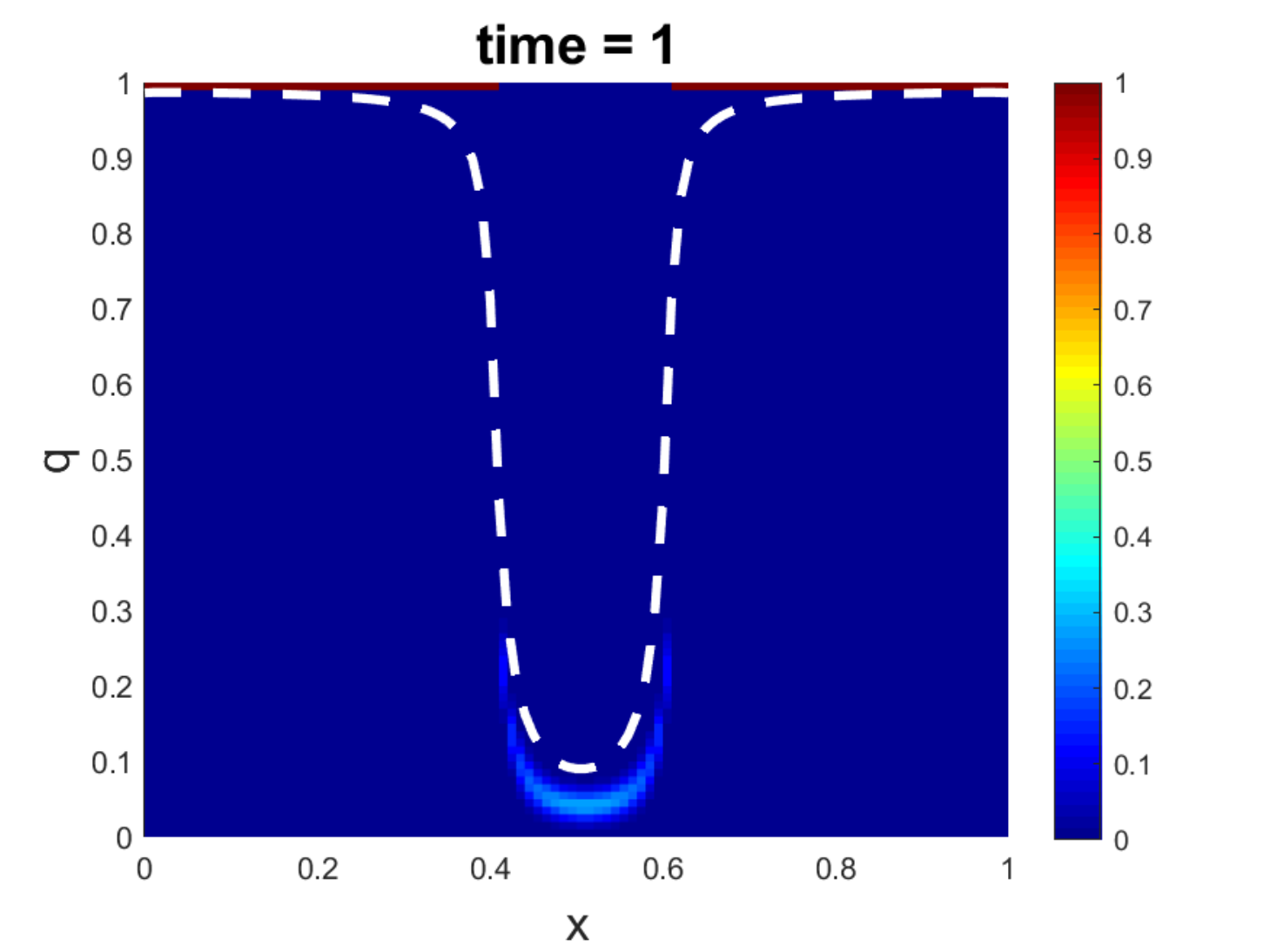}
\put(-22,45){IC1}
\put(-30,35){$\gamma = 100$}
\end{overpic}
\begin{overpic}[width=0.25\textwidth,grid=false,tics=10]{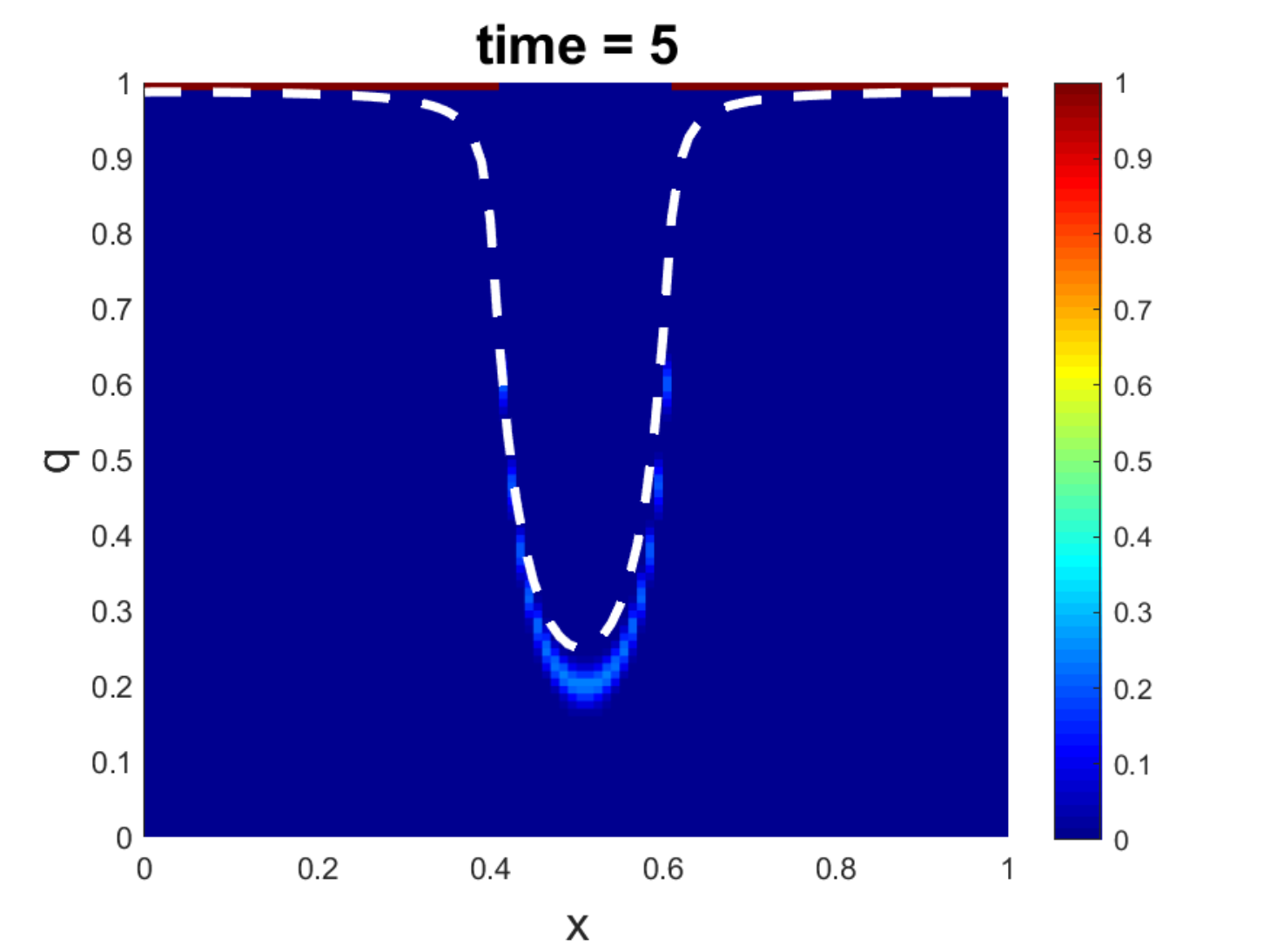}
\end{overpic} 
\begin{overpic}[width=0.25\textwidth,grid=false,tics=10]{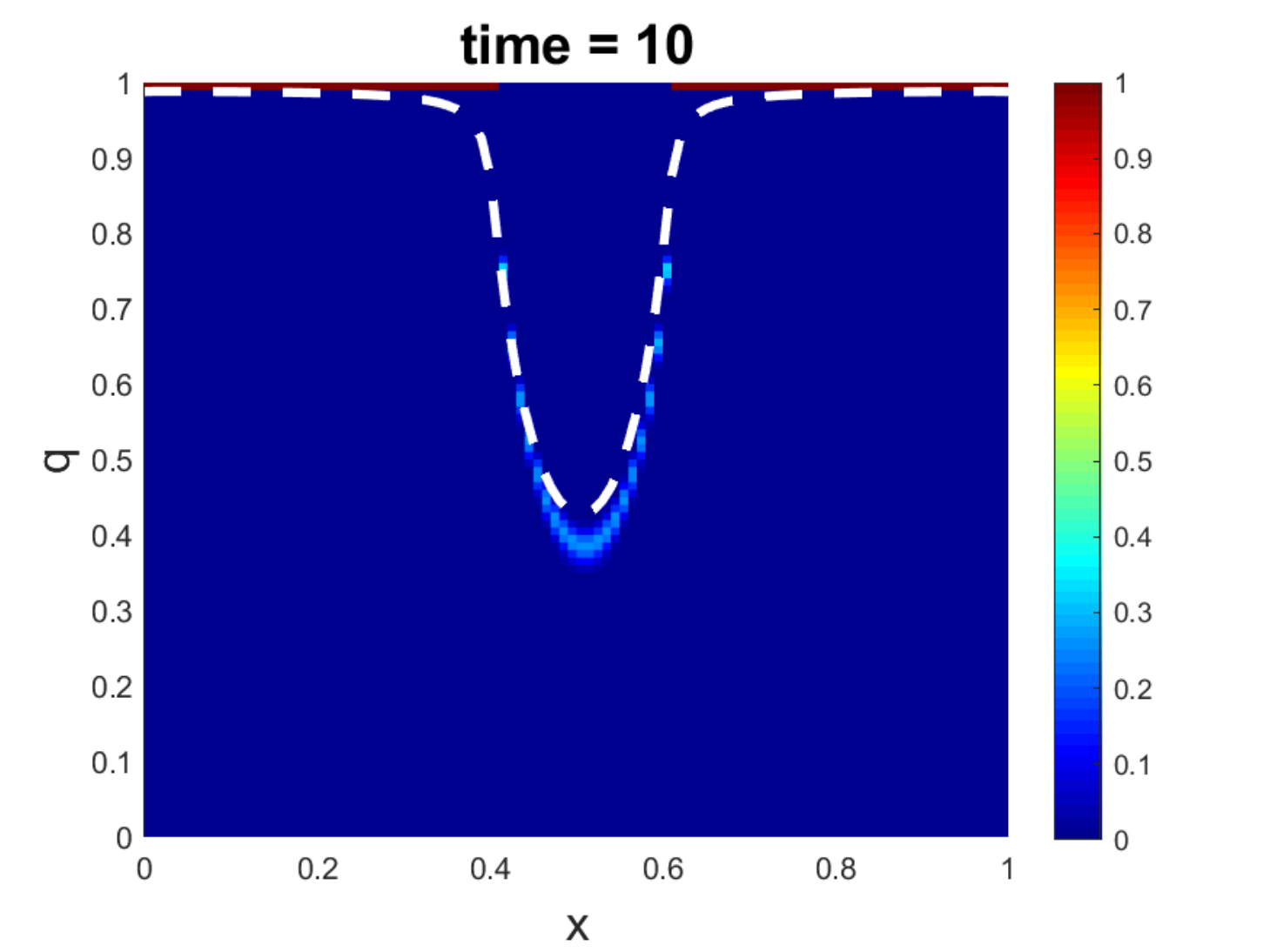}
\end{overpic} 
\begin{overpic}[width=0.25\textwidth,grid=false,tics=10]{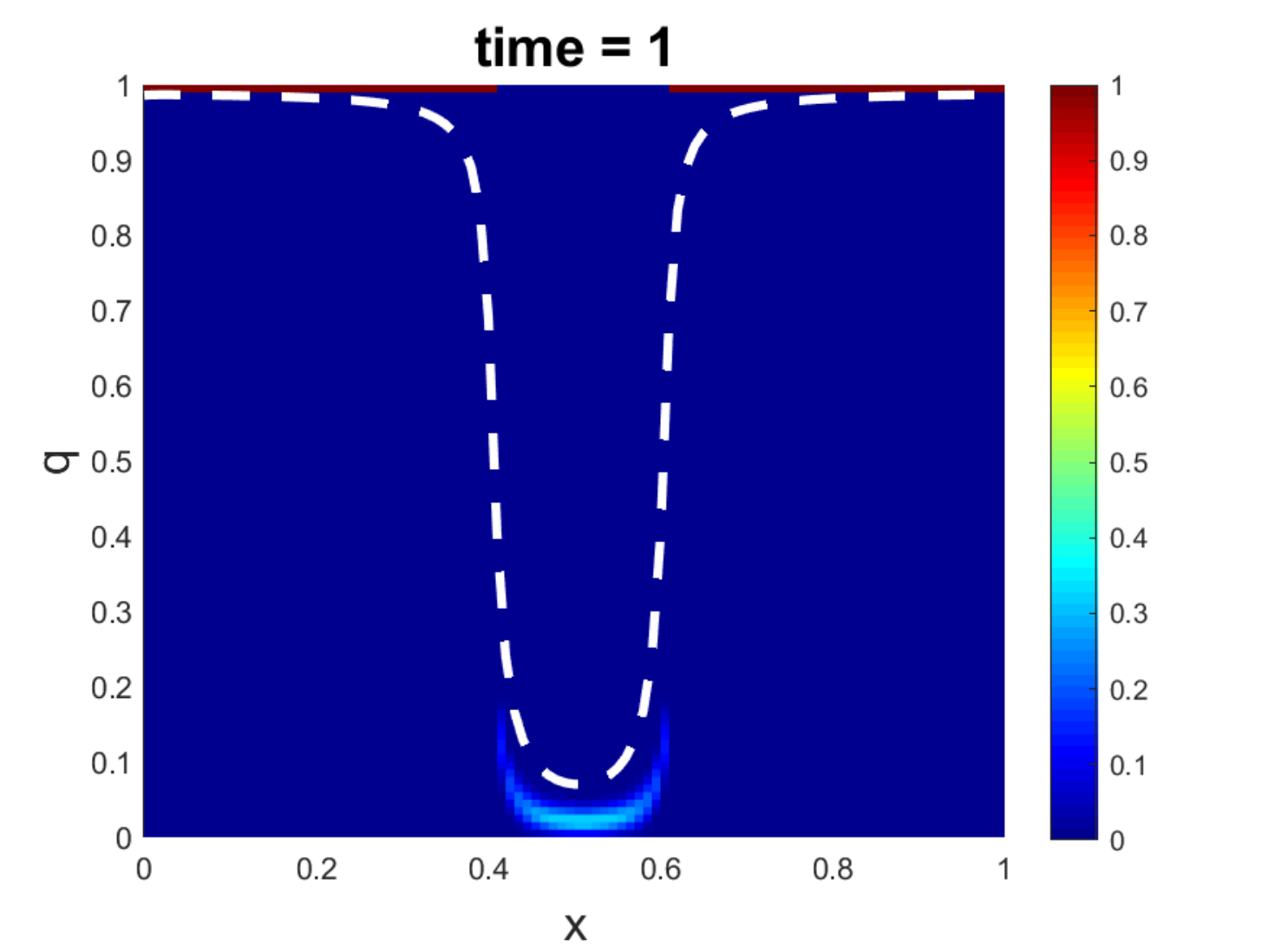}
\put(-25,45){IC1}
\put(-30,35){$\gamma = 50$}
\end{overpic}
\begin{overpic}[width=0.25\textwidth,grid=false,tics=10]{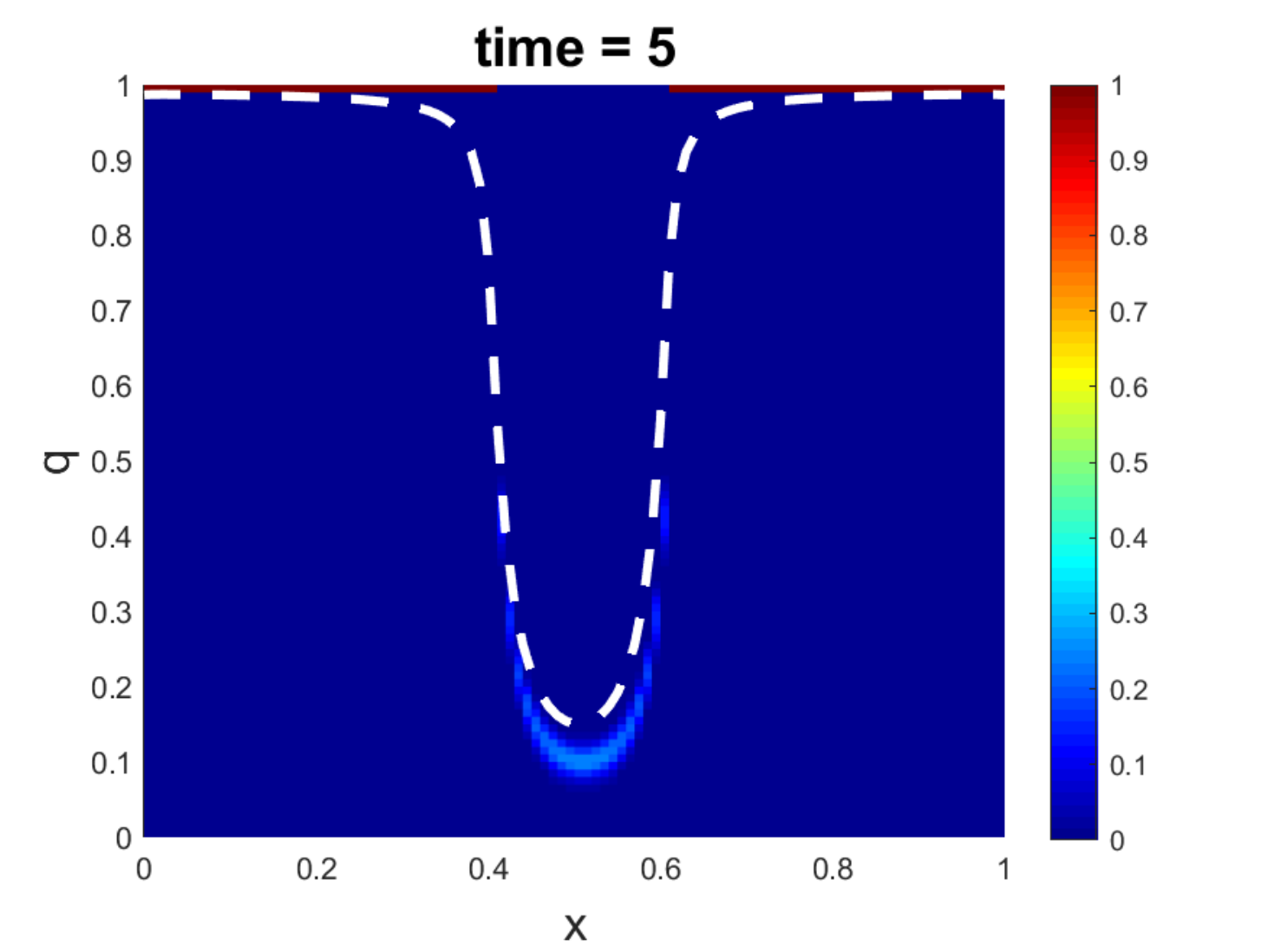}
\end{overpic} 
\begin{overpic}[width=0.25\textwidth,grid=false,tics=10]{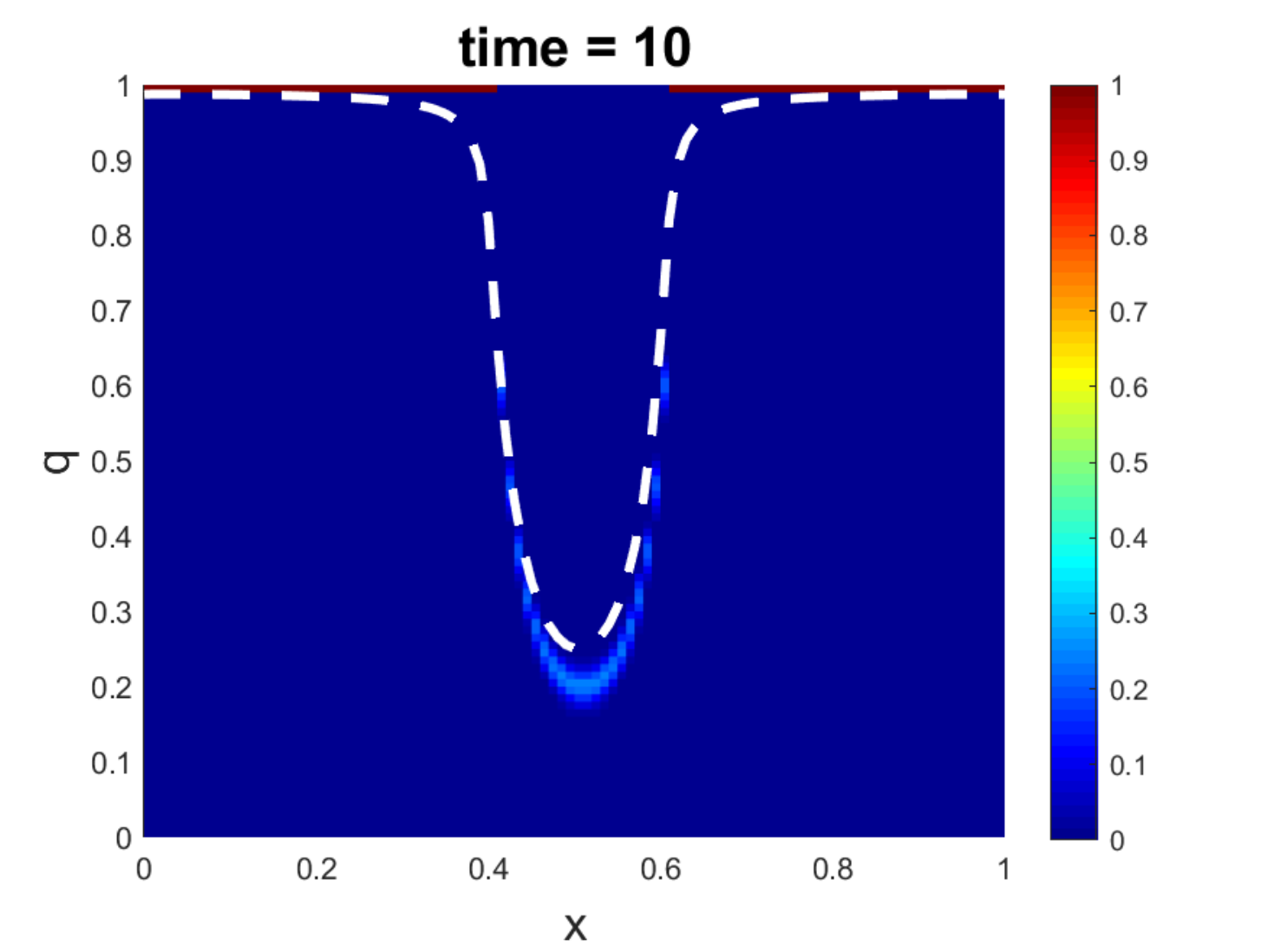}
\end{overpic} 
%\begin{overpic}[width=0.25\textwidth,grid=false,tics=10]{./possible_test2_gamma100/contagion_time_60000.pdf}
%\end{overpic}
%\begin{overpic}[width=0.25\textwidth,grid=false,tics=10]{./possible_test2_gamma100/contagion_time_300000.pdf}
%\end{overpic}
%\begin{overpic}[width=0.25\textwidth,grid=false,tics=10]{./possible_test2_gamma100/contagion_time_600000.pdf}
%\end{overpic}
\begin{overpic}[width=0.25\textwidth,grid=false,tics=10]{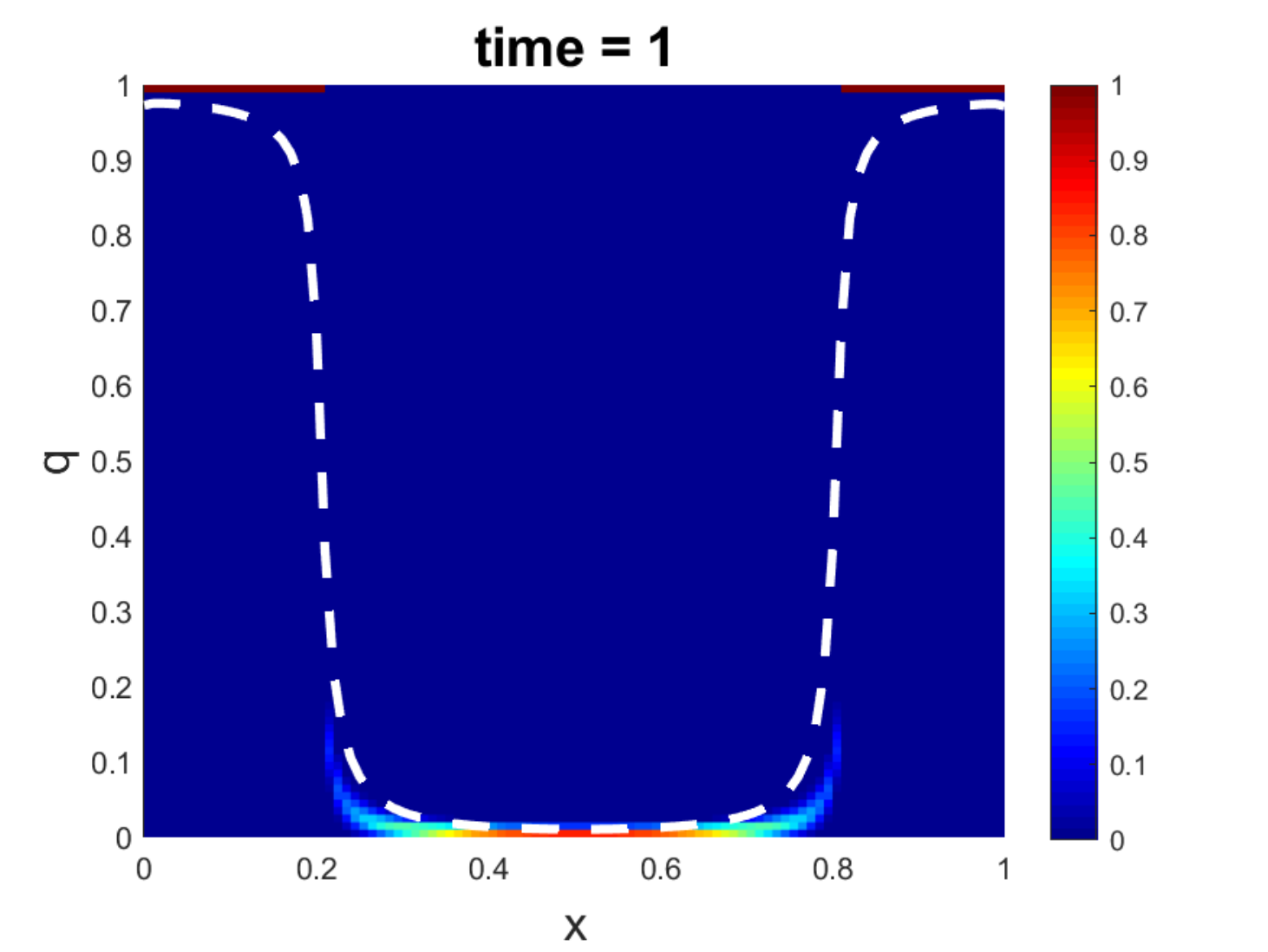}
\put(-25,45){IC2}
\put(-30,35){$\gamma = 50$}
\end{overpic}
\begin{overpic}[width=0.25\textwidth,grid=false,tics=10]{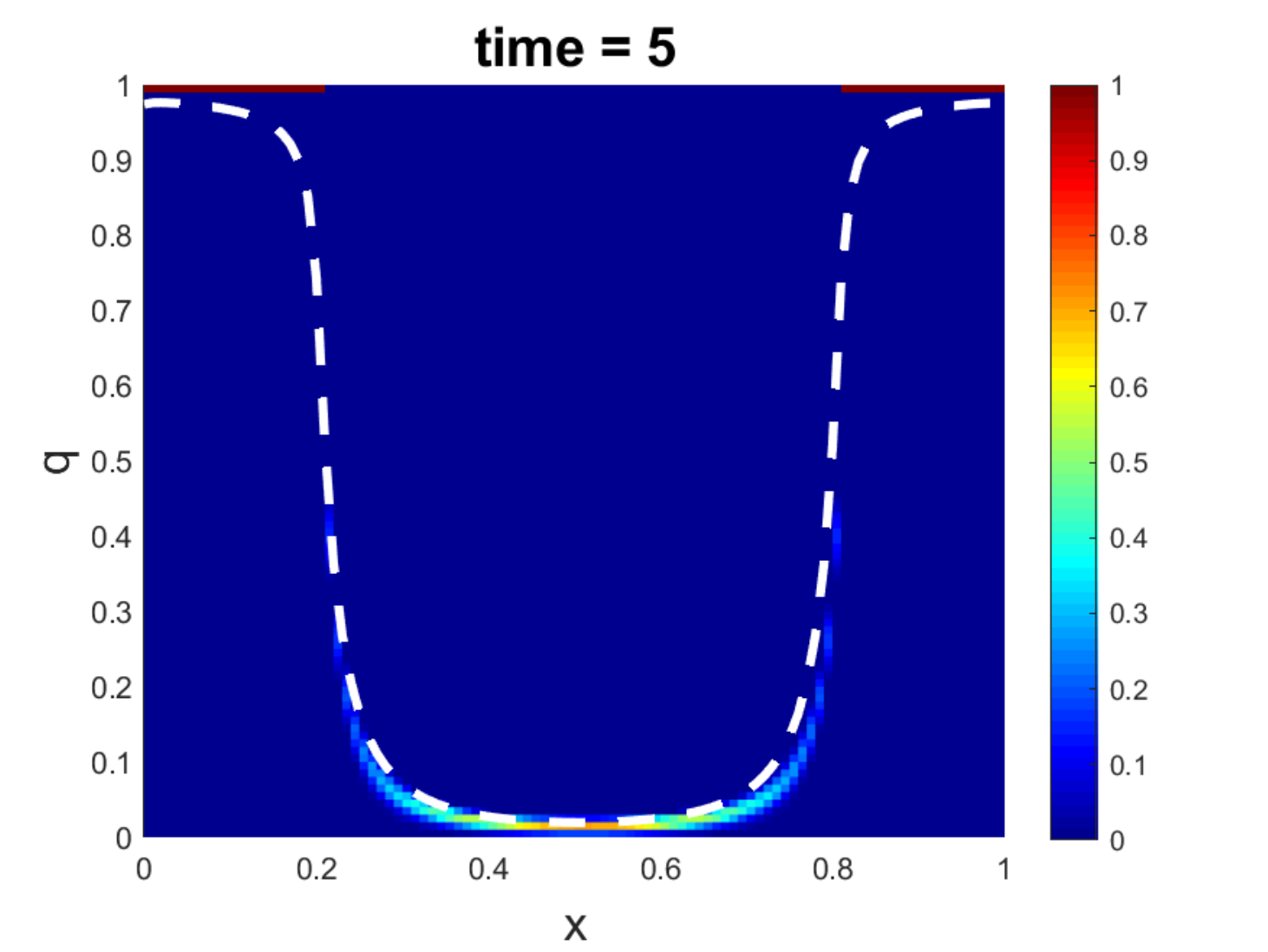}
\end{overpic}
\begin{overpic}[width=0.25\textwidth,grid=false,tics=10]{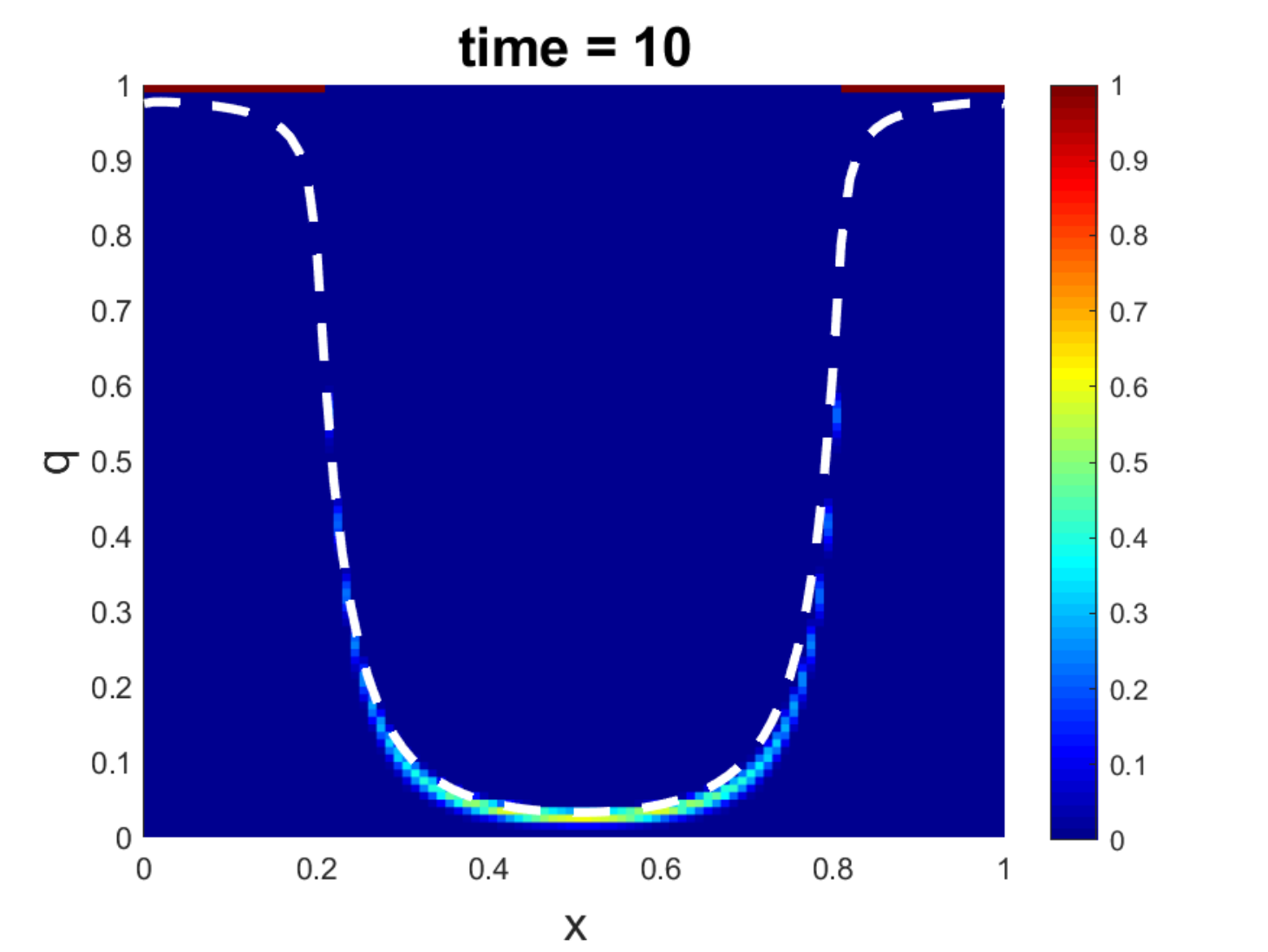}
\end{overpic}
\caption{Tests with $v = 0$: evolution of the distribution density $h$ for initial condition IC1 with $\gamma = 100$ (top)
and  $\gamma = 50$ (center), and for initial condition IC2 with $\gamma = 50$ (bottom).
The white dashed line represents $q^*$.}
\label{fig:v0}
\end{figure}

Fig.~\ref{fig:v0} shows the evolution of the distribution density $h$ for initial condition IC1 with $\gamma = 100, 50$
and for initial condition IC2 with $\gamma = 50$. We see that the level of exposure of the central group of
healthy people in IC1 increases quickly. It increases faster the closer people are to the group of sick people
and the larger $\gamma$ is. Parameter $\gamma$ plays a central role in the spreading of the disease
and would have to be carefully tuned in the future for more realistic applications. The rise in the level of exposure
is much slower for the simulation with initial condition IC2. Compare center and bottom rows in Fig.~\ref{fig:v0}.
In particular, we notice the increase in $q$ is very small for the centrally located group of healthy people, as we expected.

This first set of tests was meant to verify our implementation of method described in 
Sec.~\ref{sec:contagion_full_disc} and to check that the disease spreading term in eq.~\eqref{eq:Nkineticsystem} 
(i.e., the third term on the left-hand side) produced the expected outcomes.
Next, we are going to get people in motion. 

{\bf Tests with $\boldsymbol{v = 1}$ m/s.} We assign to all people walking direction $\theta = 0$, as if
they were headed to an exit located at $x = 10$ m. Once spreading people have left the domain, we assume they
cannot spread the disease to the people in the domain anymore. 
We consider IC1 and IC2, and set $\gamma = 50$.

\begin{figure}[htb]
\centering
\begin{overpic}[width=0.25\textwidth,grid=false,tics=10]{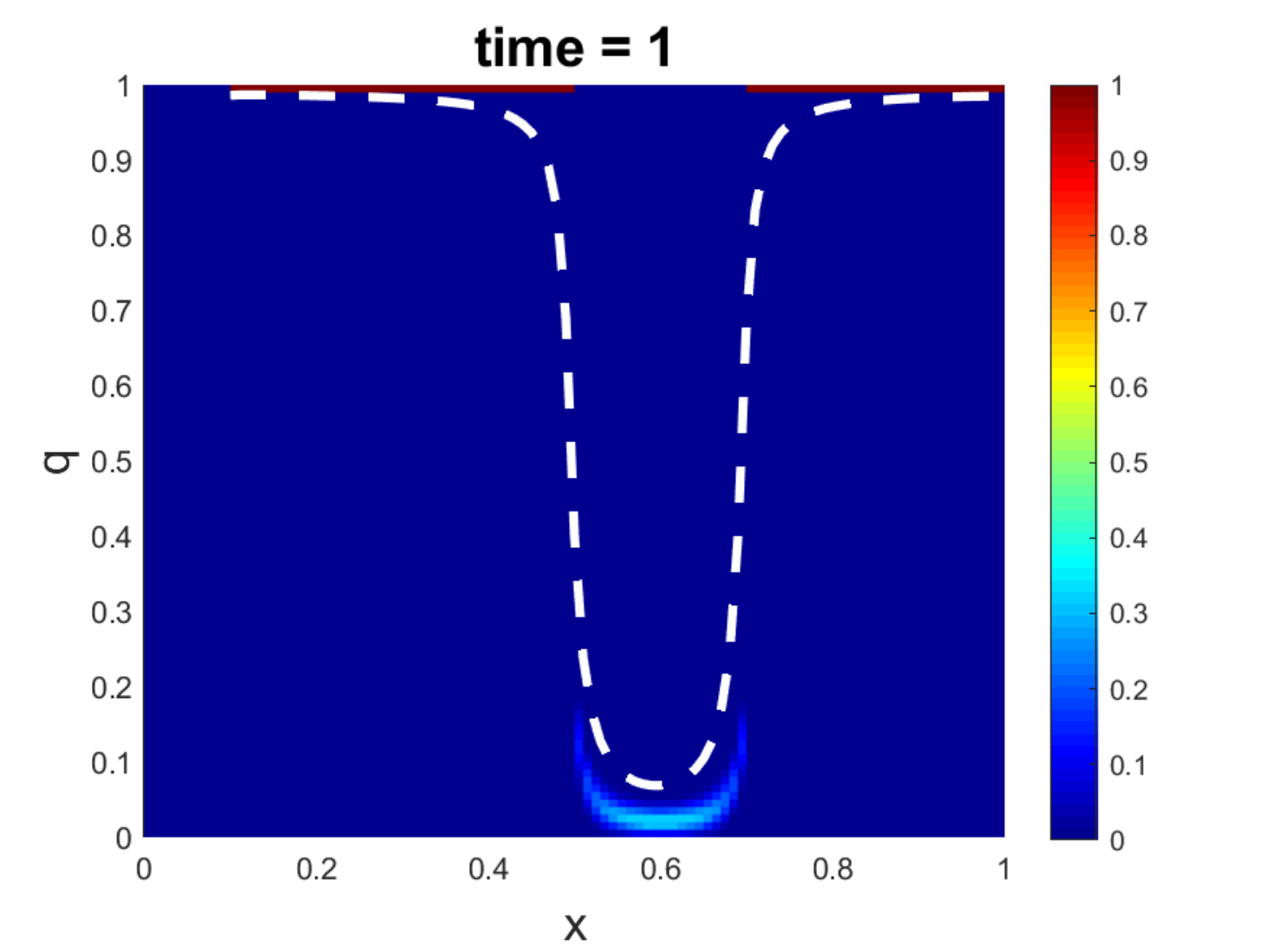}
\put(-25,37){IC1}
\end{overpic}
\begin{overpic}[width=0.25\textwidth,grid=false,tics=10]{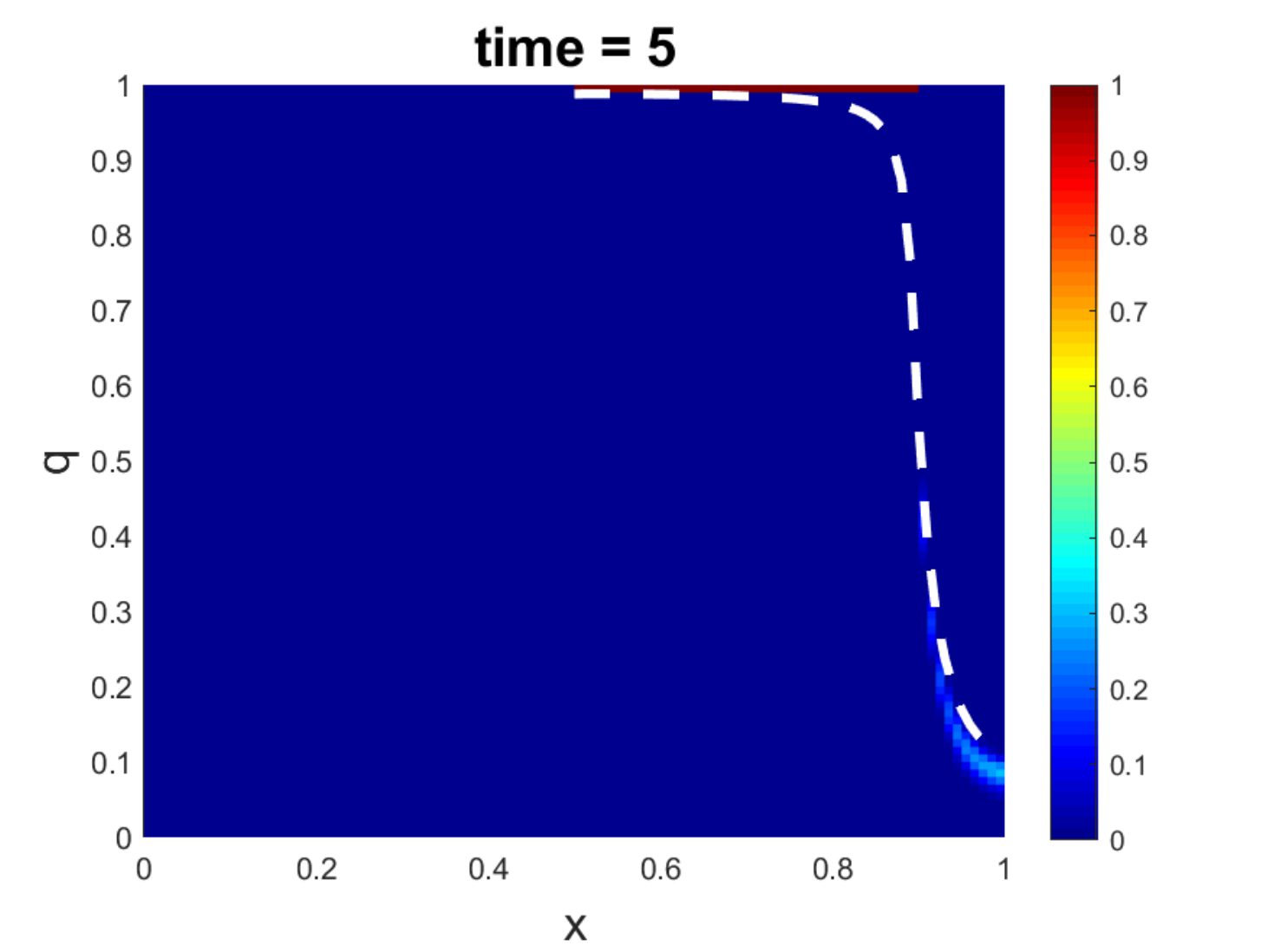}
\end{overpic} 
\begin{overpic}[width=0.25\textwidth,grid=false,tics=10]{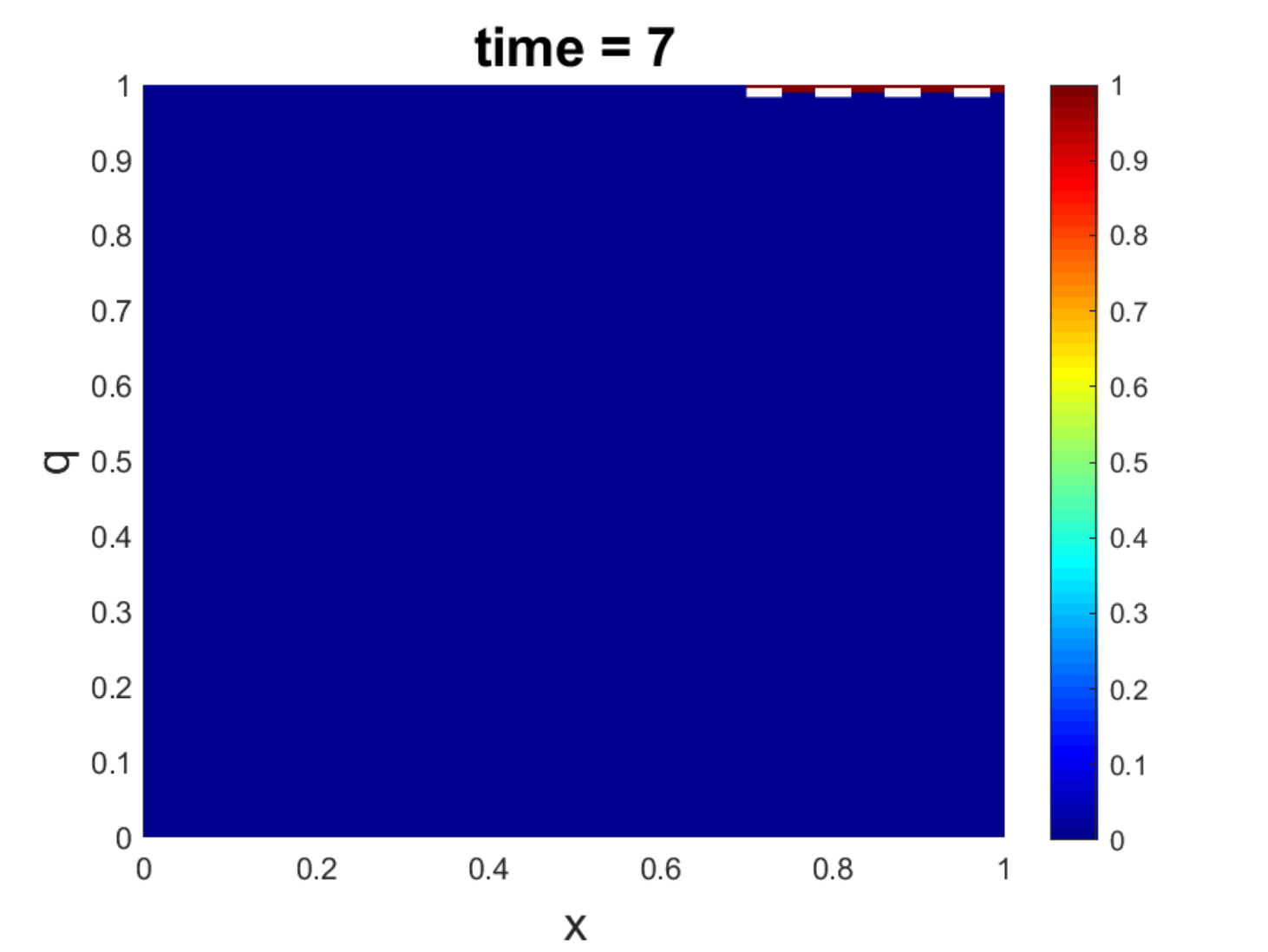}
\end{overpic} 
\begin{overpic}[width=0.25\textwidth,grid=false,tics=10]{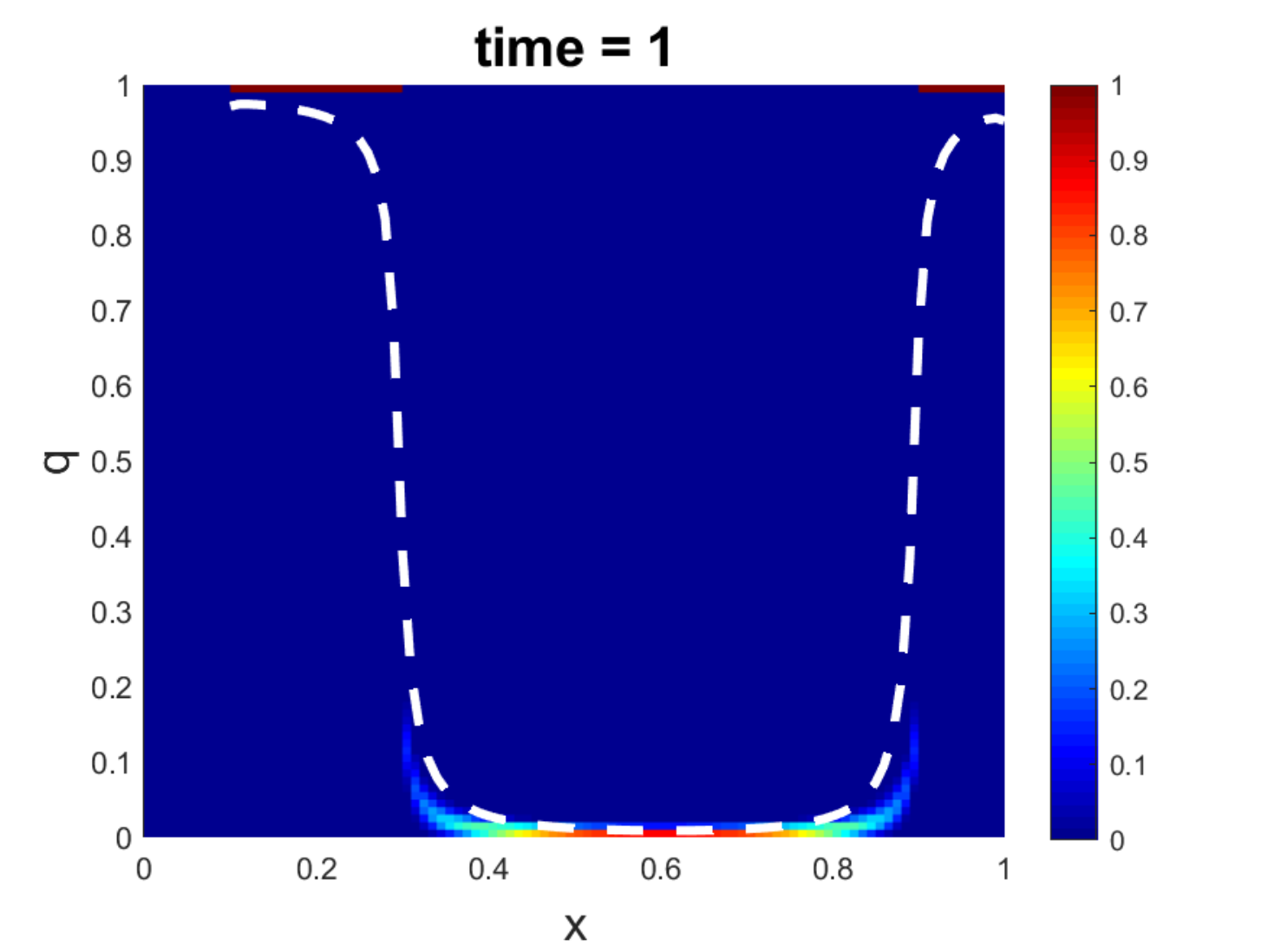}
\put(-25,37){IC2}
\end{overpic}
\begin{overpic}[width=0.25\textwidth,grid=false,tics=10]{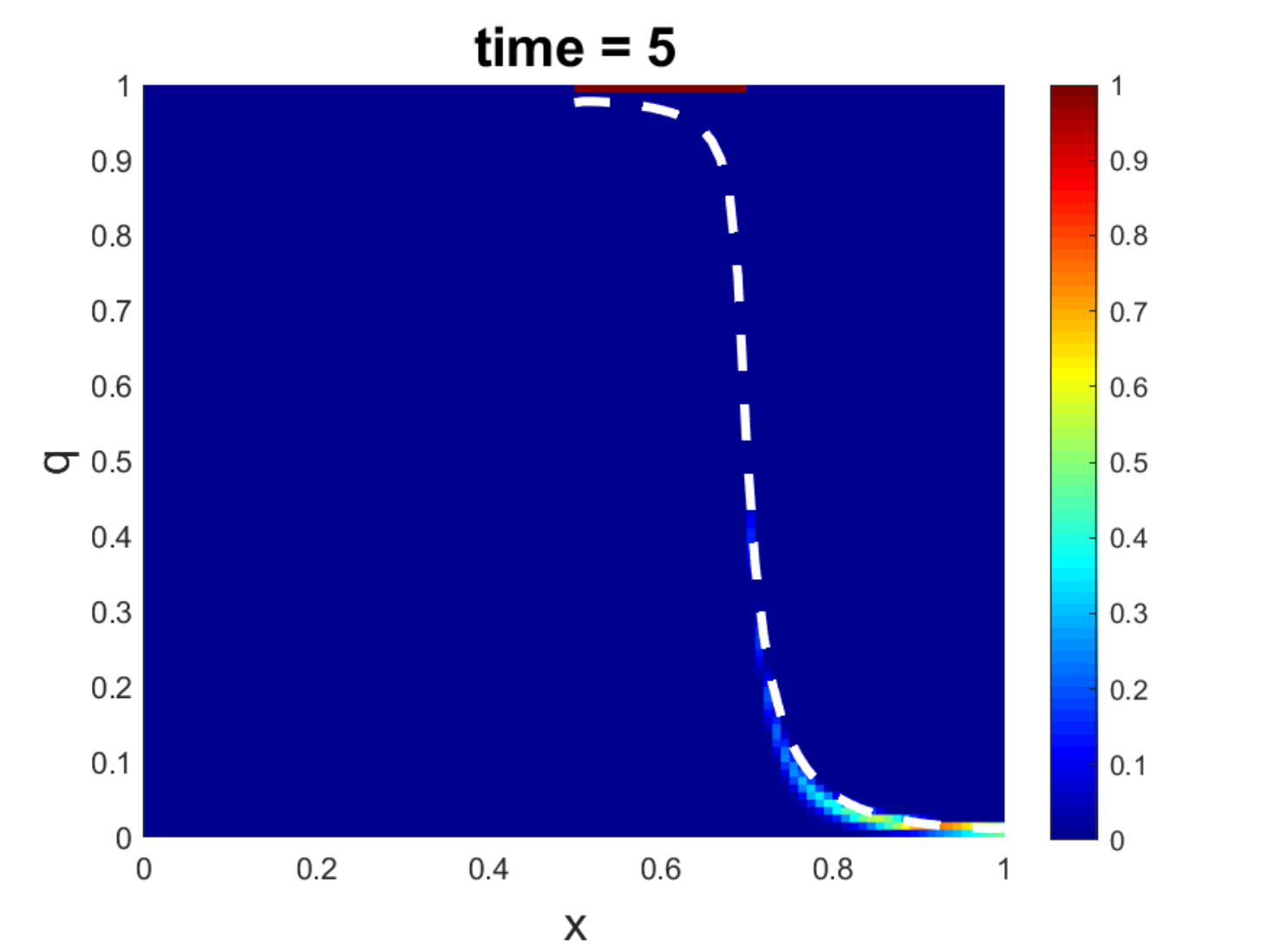}
\end{overpic} 
\begin{overpic}[width=0.25\textwidth,grid=false,tics=10]{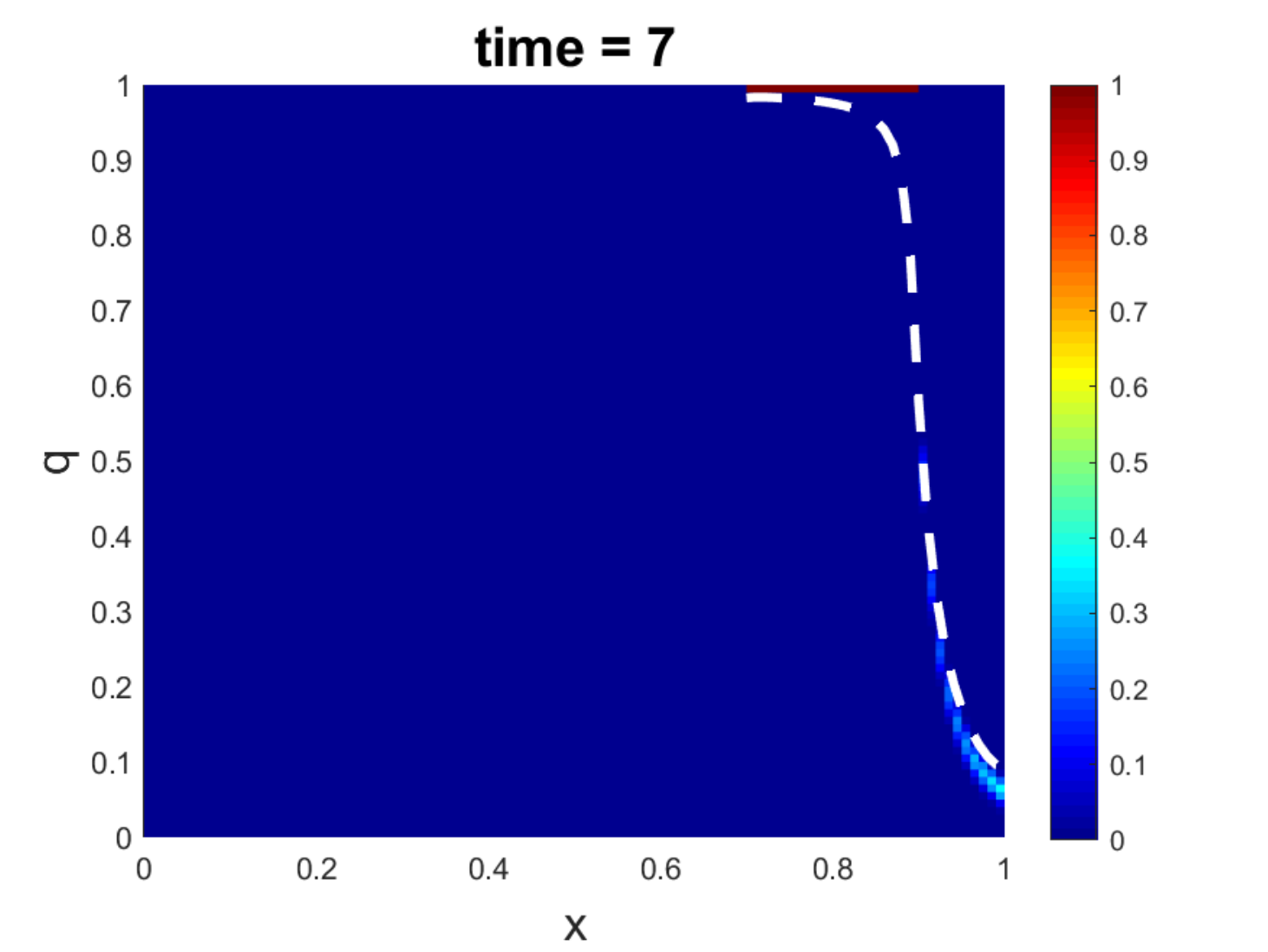}
\end{overpic} 
\caption{Tests with $v = 1$ m/s: evolution of the distribution density $h$ for initial condition IC1 (top) and IC2 (bottom). In both cases, we set $\gamma = 50$. The white dashed line represents $q^*$.} \label{fig:v1_IC}
\end{figure}

Fig.~\ref{fig:v1_IC} shows the evolution of the distribution density $h$ for initial conditions IC1 and
IC2. We observe that the motion contributes to lowering the exposure level in the both cases, 
since some of the spreading people leave the domain first. Compare the top and bottom row of Fig.~\ref{fig:v1_IC} 
with the central and bottom row of Fig.~\ref{fig:v0}. 

\begin{figure}[ht!]
\centering
\begin{overpic}[width=0.32\textwidth,grid=false,tics=10]{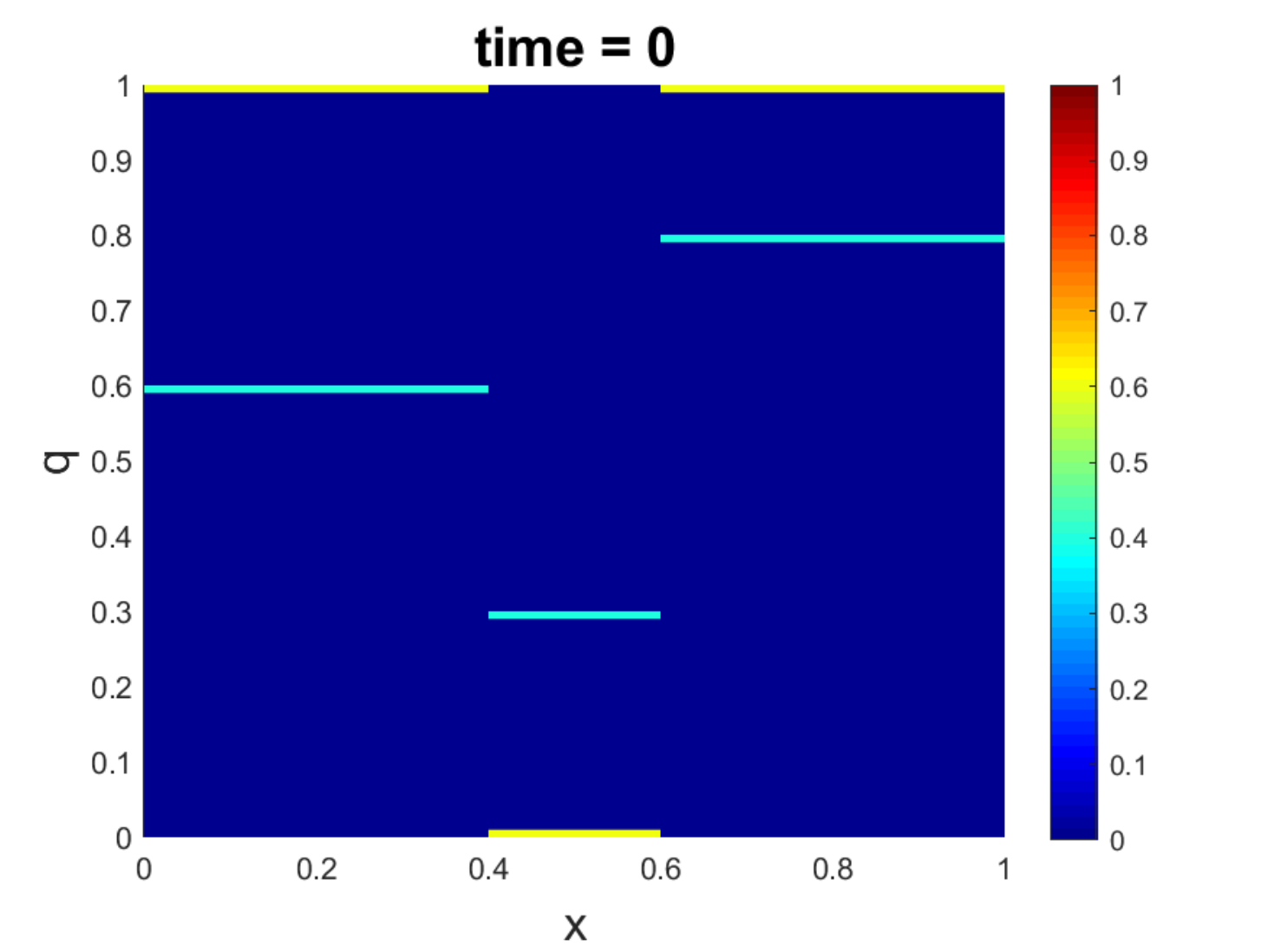}
\end{overpic} 
\begin{overpic}[width=0.32\textwidth,grid=false,tics=10]{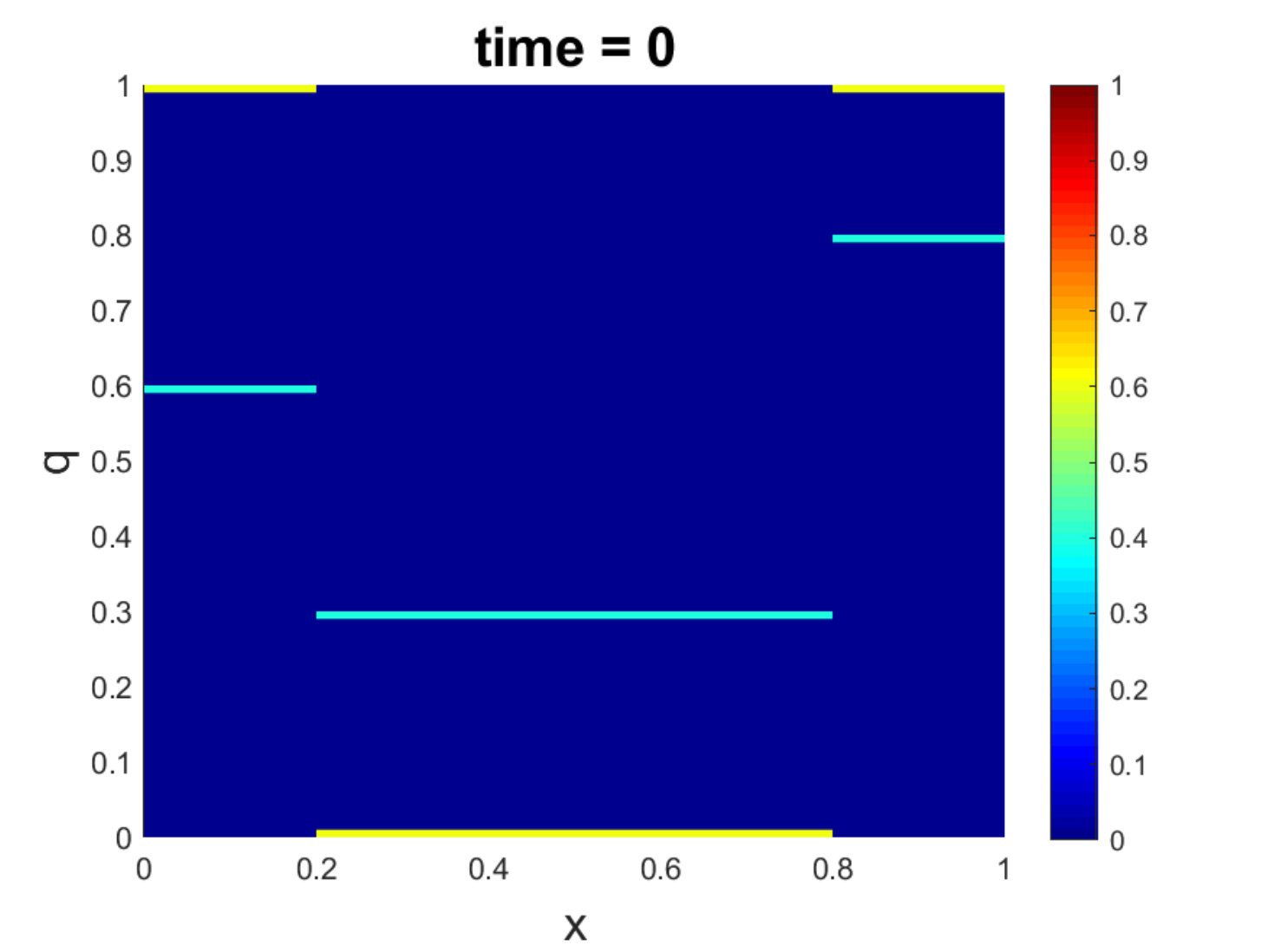}
\end{overpic} 
\caption{Tests with $v = 1$ m/s: initial conditions IC1-bis (left) and IC2-bis (right).}\label{fig:ICbis}
\end{figure}

\begin{figure}[htb]
\centering
\begin{overpic}[width=0.25\textwidth,grid=false,tics=10]{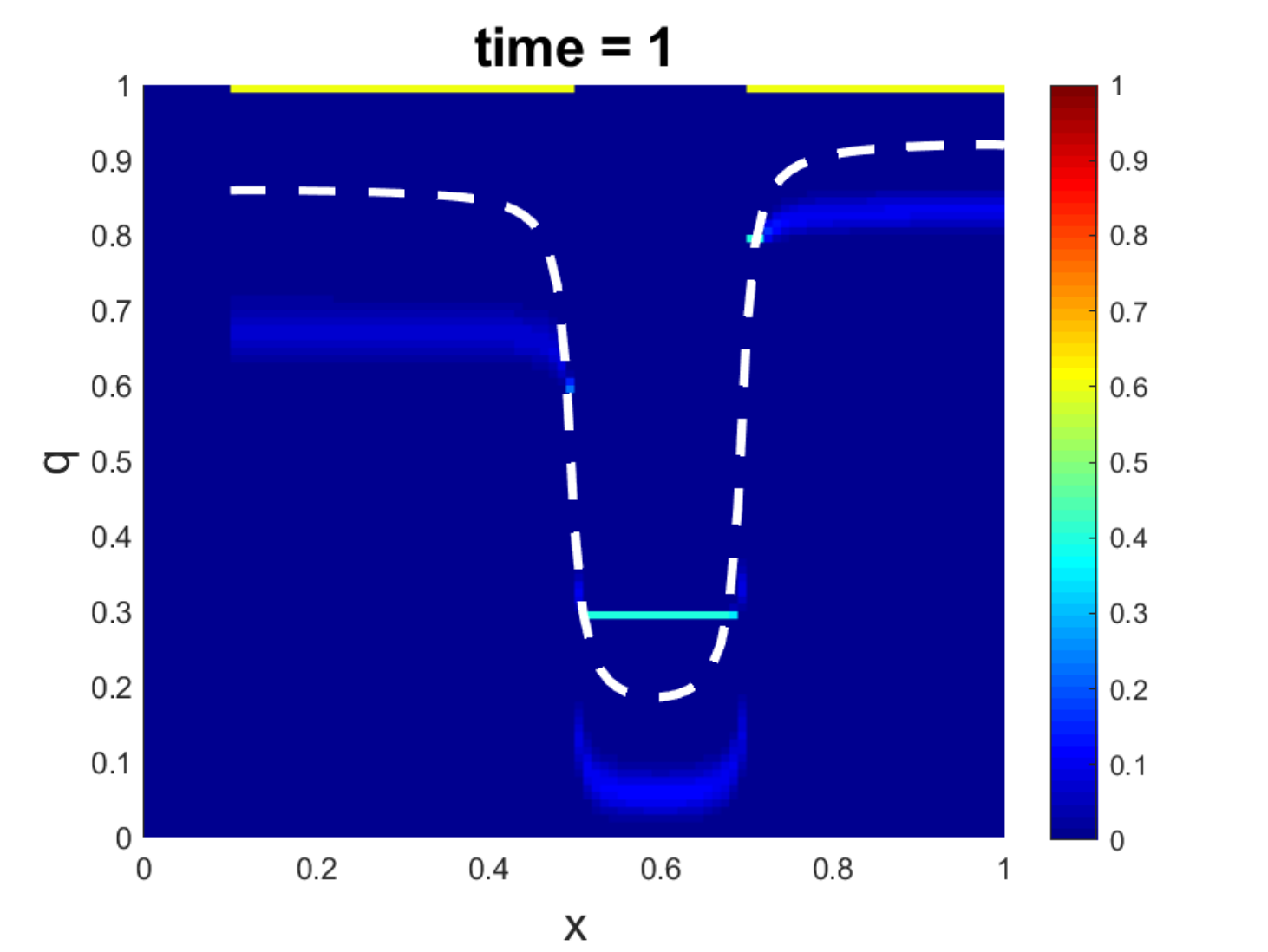}
\put(-30,37){IC1-bis}
\end{overpic}
\begin{overpic}[width=0.25\textwidth,grid=false,tics=10]{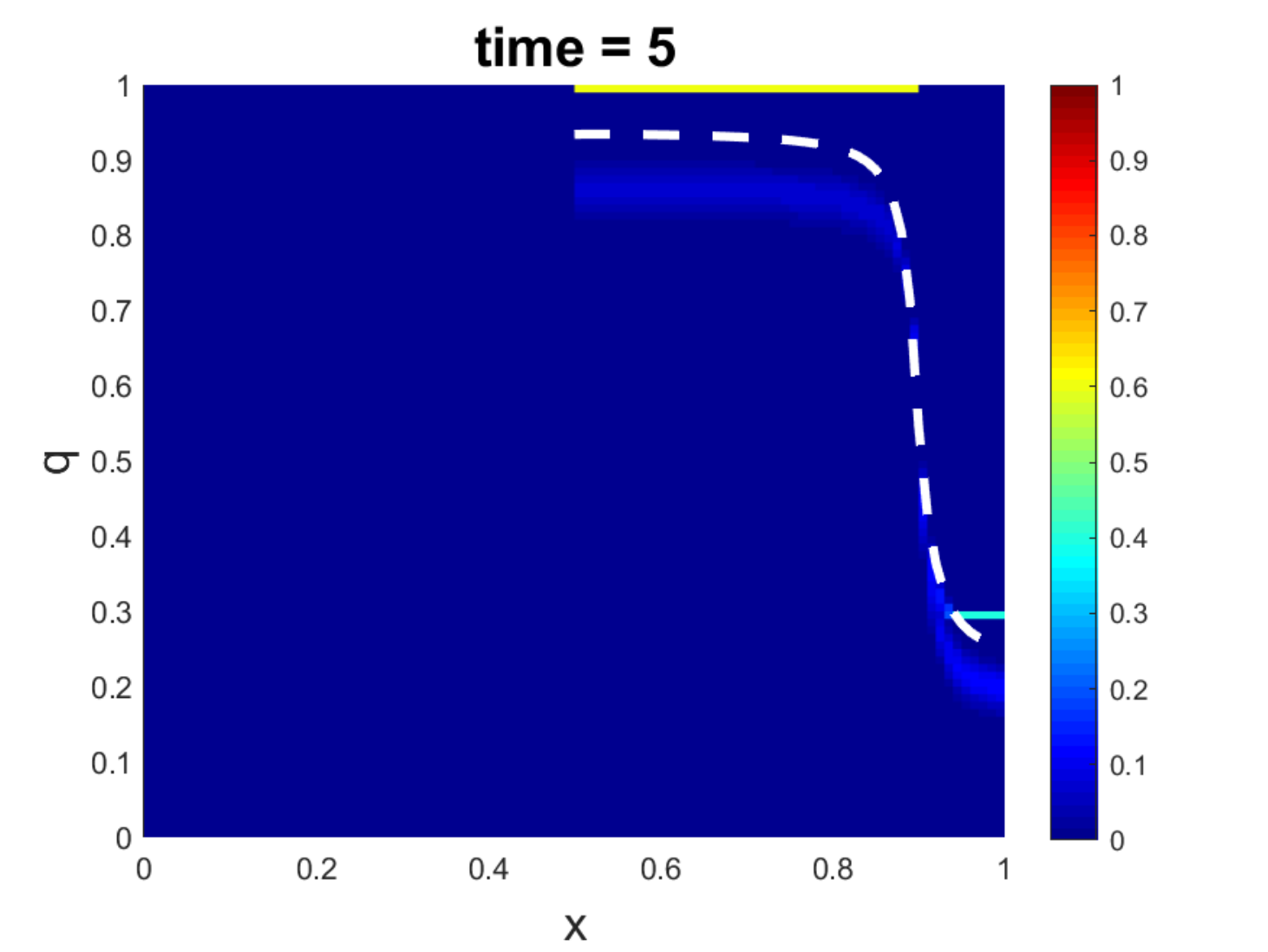}
\end{overpic} 
\begin{overpic}[width=0.25\textwidth,grid=false,tics=10]{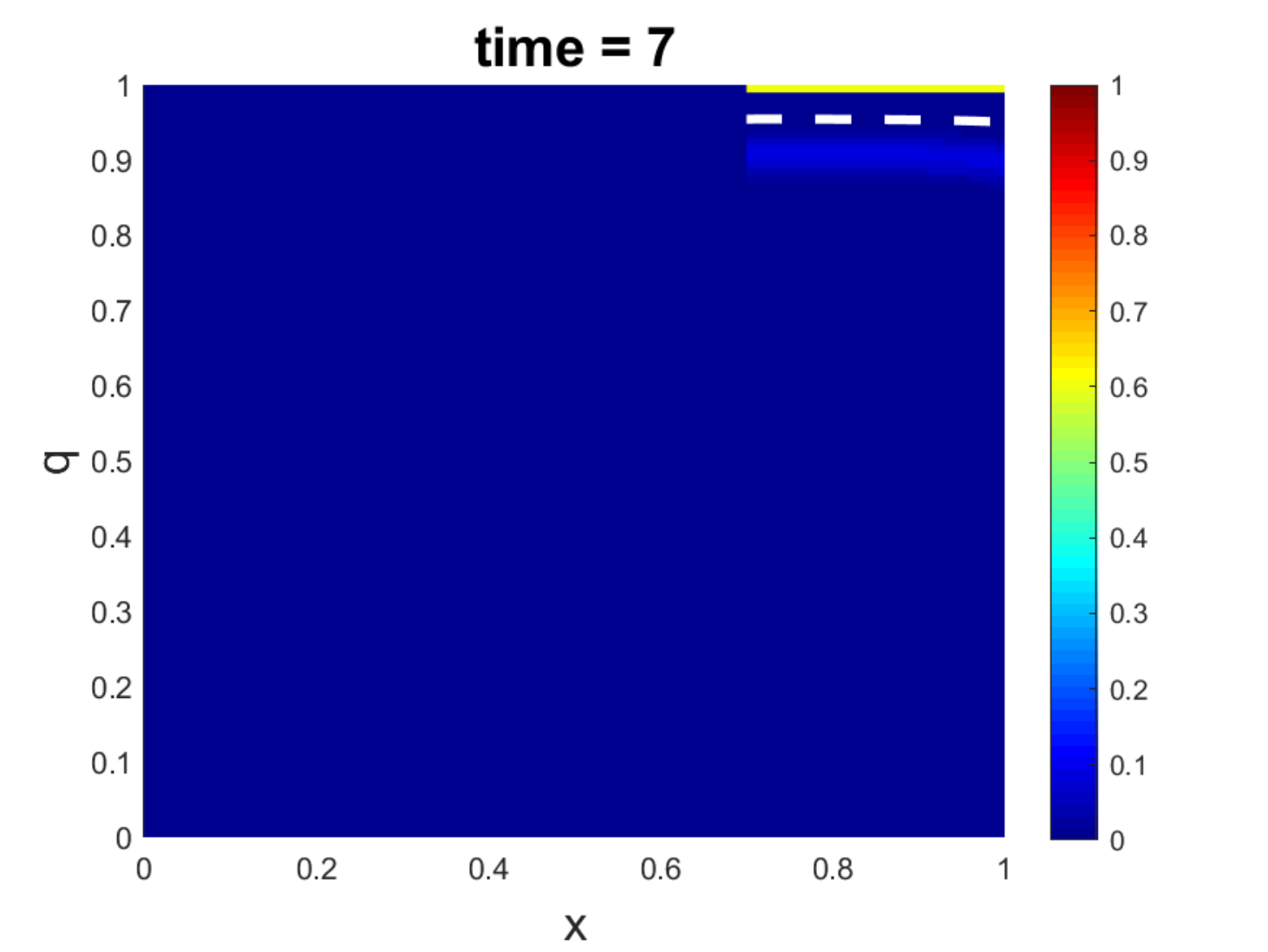}
\end{overpic} 
 \begin{overpic}[width=0.25\textwidth,grid=false,tics=10]{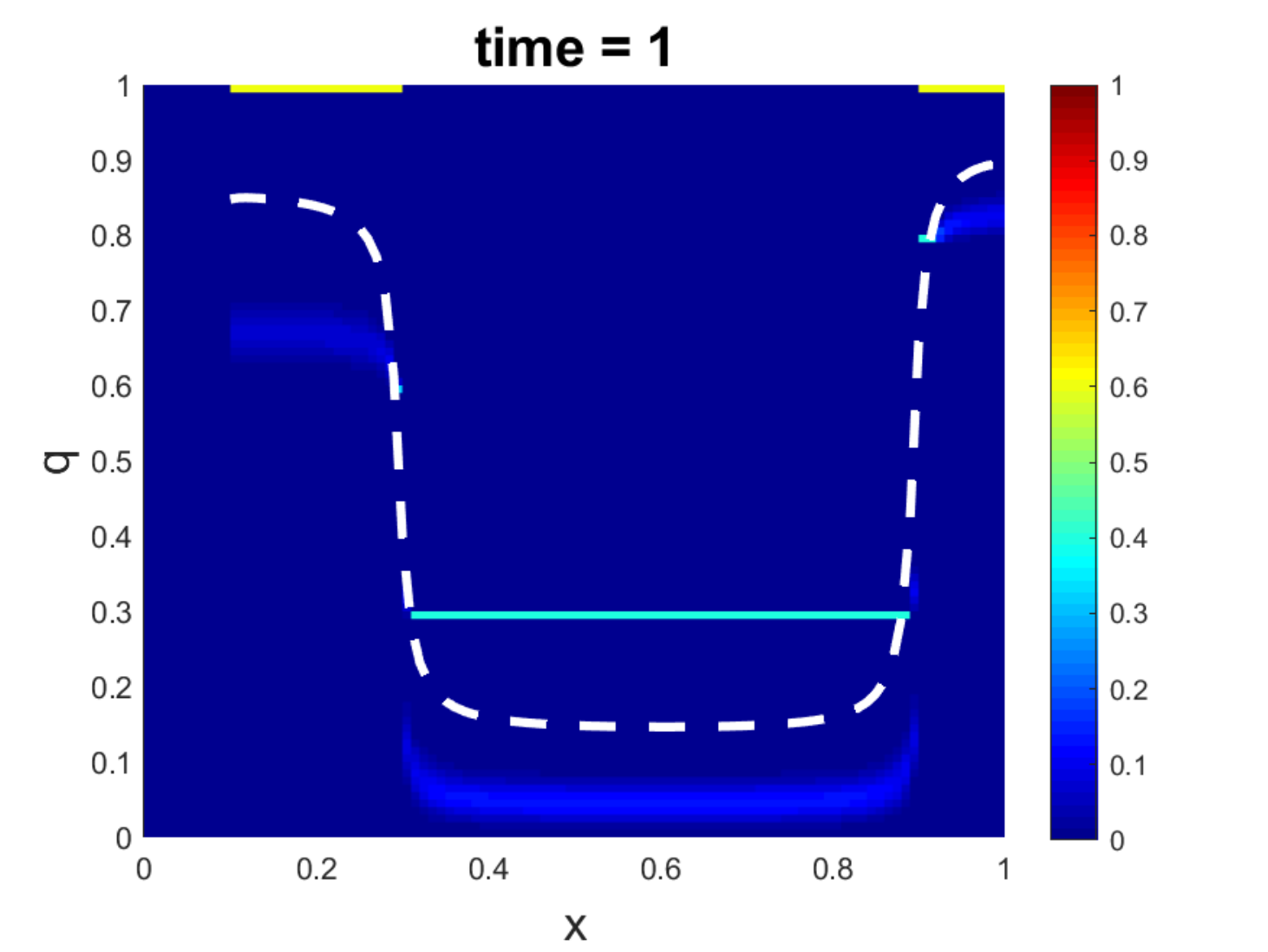}
 \put(-30,37){IC2-bis}
\end{overpic}
\begin{overpic}[width=0.25\textwidth,grid=false,tics=10]{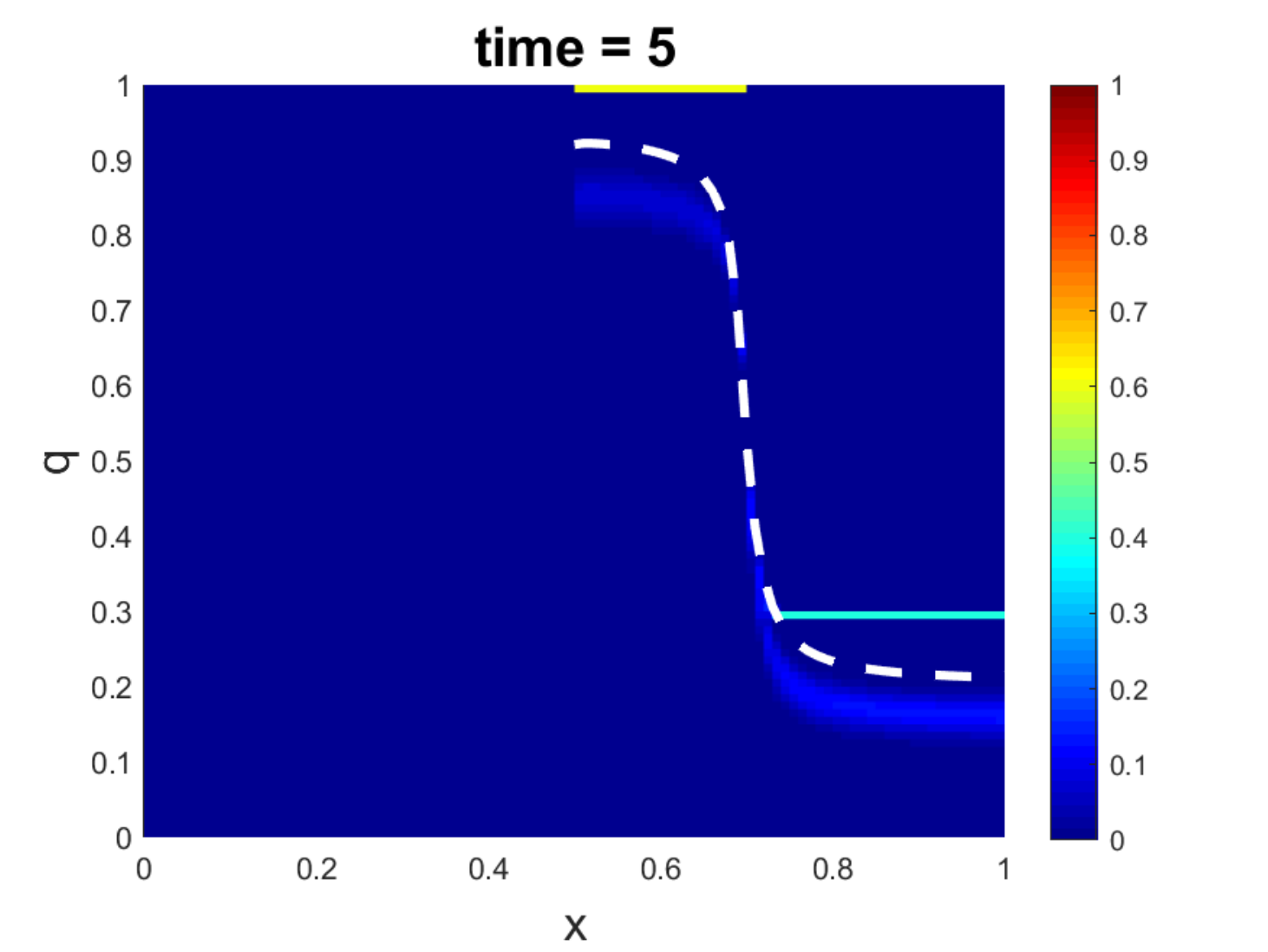}
\end{overpic} 
\begin{overpic}[width=0.25\textwidth,grid=false,tics=10]{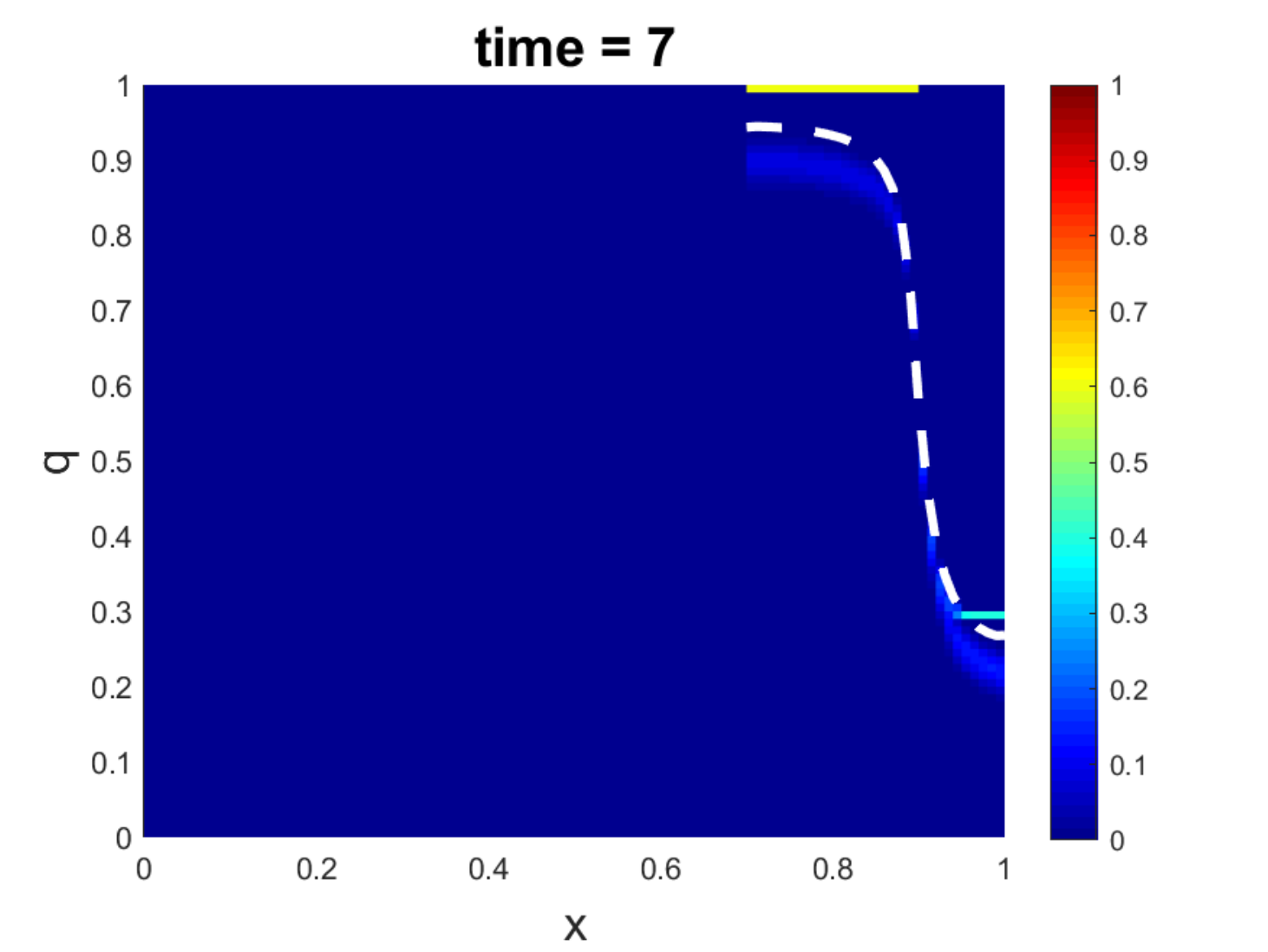}
\end{overpic} 
\caption{Tests with $v = 1$ m/s: evolution of the distribution density $h$ for initial condition IC1-bis (top)
and IC2-bis (bottom). In both cases, we set $\gamma = 50$. The white dashed line represents $q^*$.} \label{fig:v1_ICbis}
\end{figure}

Finally, we experiment with a slight modification of the initial conditions to show that
our model can handle scenarios with uncertainty. The initial conditions
are changed:
\begin{itemize}
\item[-] IC1-bis: people are positioned like in IC1 but the probabilities of finding people with $q = 1$
and $q = 0$ is reduced from 100\% to 60\% and another value of $q$ for a given $x$ is assigned. See Fig.~\ref{fig:ICbis} (left panel).
\item[-] IC2-bis: people are positioned like in IC2 but the probabilities of finding people with $q = 1$
and $q = 0$ is reduced from 100\% to 60\% and another value of $q$ for a given $x$ is assigned. See Fig.~\ref{fig:ICbis} (right panel).
\end{itemize}
Fig.~\ref{fig:v1_ICbis} shows the evolution of the distribution density $h$ for initial conditions IC1-bis and
IC2-bis. 

%% conclusion  or discussion %%%%%%%%%%%%%%%%%%%%%
\section{Conclusion}\label{sec:concl}
%%%%%%%%%%%%%%%%%%%%%%%%%%%%
This paper is divided into two parts. In the first part we presented a kinetic type model for crowd dynamics, 
while in the second part we introduced a simplified model for disease contagion in a crowd
walking through a confined environment. 

Kinetic (or mesoscopic) approaches to simulate the motion of medium-sized crowds are appealing 
because of their flexibility in accounting for multiple interactions (hard to achieve in microscopic models) 
and heterogeneous behavior in people (hard to achieve in macroscopic models). The particular kinetic model
we chose was also shown to compare favorably with experimental data for a medium-sized population.
Previously, this model had been used to simulate simple scenarios such as evacuation with a room. 
In this paper, we showed that realistic scenarios, such as passengers walking in an airport terminal, 
can be handled as well. 

The simplifying assumptions that we used in the model for disease contagion is that
people's walking speed and direction are given. The disease spreading is modeled
using three main ingredients: an additional variable that denotes the level of exposure to people 
spreading the disease, a parameter that describes the contagion interaction strength, 
and a kernel function that is a decreasing function of the distance between a person and a spreading individual.
We tested the proposed contagion model and numerical approach on simple 1D problems.

The obvious next step is to combine the kinetic type model for crowd dynamics in the first part
of the paper with the disease contagion model in order to drop the simplifying assumption, i.e.~walking 
speed and direction are provided by the model instead of being given.

\section*{Acknowledgements}
This work has been partially supported by NSF through grant DMS-1620384.

%\clearpage
\bibliographystyle{plain}
\bibliography{contagion_model}

\end{document}